\documentclass[11pt]{article}
\usepackage{amsmath}
\usepackage{amssymb}
\usepackage{amsthm}
\usepackage{url}
\usepackage{color}
\usepackage[pstarrows]{pict2e}
\numberwithin{equation}{section}
\newtheorem{theorem}{Theorem}[section]
\newtheorem{proposition}[theorem]{Proposition}
\newtheorem{Remark}[theorem]{Remark}
\newenvironment{remark}{\begin{Remark}\rm}{\end{Remark}}
\newcommand\Proof{\noindent{\bf Proof}\quad}
\newcommand\nonu{\nonumber}
\newcommand\sLP{\\[\smallskipamount]}
\newcommand\mLP{\\[\medskipamount]}
\newcommand\bLP{\\[\bigskipamount]}
\newcommand\RR{\mathbb{R}}
\newcommand\al\alpha
\newcommand\be\beta
\newcommand\ga\gamma
\newcommand\de\delta
\newcommand\la\lambda
\newcommand\si\sigma
\newcommand\iy\infty
\newcommand{\hyp}[5]{\,\mbox{}_{#1}F_{#2}\!\left(
  \genfrac{}{}{0pt}{}{#3}{#4};#5\right)}
\newcommand\thalf{\tfrac12}
\newcommand\wrt{with respect to }

\begin{document}
\title{The Askey scheme as a four-manifold with corners}
\author{Tom H. Koornwinder\bLP
{\small  Korteweg-de Vries Institute, University of Amsterdam,}\\
{\small{\tt T.H.Koornwinder@uva.nl}}}
\date{}
\maketitle
\begin{abstract}
Racah and Wilson polynomials with dilated and translated argument are
reparametrized such that the polynomials are continuous in the
parameters
as long as these are nonnegative, and such that restriction of one or more
of the new parameters to zero yields orthogonal polynomials lower in the
Askey scheme. Geometrically this will be described as a manifold with
corners.
\end{abstract}
%
%
\section{Introduction}
A graphical scheme describing the families of
hypergeometric orthogonal polynomials and the limit relations between them
was published by Askey \& Wilson in 1985 as an Appendix to
their famous memoir \cite{14} that introduced the Askey-Wilson polynomials.
Immediately afterwards, the scheme was slightly extended
for the continuous Hahn polynomials, see
\cite{15}, \cite{16}, and it remained stable since then.
The scheme was named after Askey by J. Labelle \cite{17}
(tableau d'Askey or Askey scheme), who also made a big wall poster
displaying the scheme. As often told by Askey, he got the idea for
this scheme already at an Oberwolfach meeting in 1977 on
``Combinatorics and Special Functions''. There Michael Hoare, in connection
with his lecture, distributed copies of a sheet (an extension
of \cite[p.285, Figure 2]{18}) which
contained in graphical way a part of the present Askey scheme,
and which was received very enthusiastically by the audience.

Since its introduction, the Askey scheme (given here in Figure \ref{fig:1})
has been very influential both
by its compact way of presenting the various families of
hypergeometric orthogonal polynomials and by the structuring effect
of the arrows which indicate the limit relations. Furthermore, the report
\cite{7} by Koekoek and Swarttouw has been very helpful in giving the
most important formulas for each family and each limit relation in
the scheme (and in the $q$-Askey scheme as well).

Personally, I was in particular intrigued by the arrows. In 1993 I
wrote a report \cite{19} in which, for a part of the Askey scheme,
limit relations are described together in a uniform way.
In different notation, these results are given in the present paper
in \S\ref{35} and the corresponding part of the Askey scheme is in the
first graph
of Figure \ref{fig:2}. The idea here is to reparametrize the
four-parameter family of Racah polynomials, together with a
parameter dependent translation and dilation of the argument, such that
the resulting orthogonal polynomials become continuous up to the
boundary in their dependence on the four nonnegative parameters.
Restriction to some part of the boundary yields families lying below
the Racah polynomials in the first graph of Figure \ref{fig:2}, see also
Figure \ref{fig:4}. The method of proof is by examining continuous
dependence on the parameters in the coefficients for the
three-term recurrence relation of the monic orthogonal polynomials.

The present paper extends the results of \cite{19} such that the whole
Askey scheme is covered. It is helpful to view this in a geometric way
as several four-manifolds with coordinates coming from the
parameters for the Racah or Wilson polynomials. Here we have
one four-manifold for the Racah polynomials and two distinct
four-manifolds for the Wilson polynomials. We should consider these
four-manifolds together with boundary, not just a boundary of
codimension 1, but boundary parts of all dimensions $0,1,2,3$.
Such manifolds (in arbitrary finite dimension) are known as
{\em manifolds with corners}, introduced by Cerf \cite{12}
and Douady \cite{13}. The various boundary parts will correspond
with families lying below the top level in the Askey scheme.
It turns out that for each of the two Wilson manifolds (see
Figure \ref{fig:3}) one
reparametrization is sufficient. But in order to cover the
Racah polynomials and everything below, three reparametrizations
(local charts) are needed, which partially overlap. See
Figure \ref{fig:2} for a graphical representation of what is
covered by each of the three charts.

The detailed results are quite computational. In order to obtain them,
I made heavy use of the computer algebra program
{\em Mathematica}${}^{\mbox{\footnotesize\textregistered}}$.
The Mathematica notebooks with the computations will be available on the web,
see \url{http://www.science.uva.nl/~thk/art/}.

Various papers have appeared about asymptotics in connection with the
Askey scheme, see for instance
\cite{21}, \cite{22}, \cite{4}, \cite{5}, \cite{3}, \cite{1}, \cite{2},
\cite{23}, \cite{24}.
Limits are a very special case of asymptotics. The technique
of the present paper for obtaining limits using the three-term
recurrence relation may be also relevant for further asymptotics.

Another potential application of this paper would be in the problem
to identify some explicit system of orthogonal polynomials when the
explicit coefficients in its three-term recurrence relation are given.
In joint work by the author with Swarttouw a
Maple${}^{\mbox{\footnotesize\textregistered}}$ procedure called
{\tt rec2ortho} was written which can do this job for polynomials
in the Askey scheme up to the
2-parameter level (Jacobi, Meixner, etc.), see
\url{http://www.science.uva.nl/~thk/art/software/rec2ortho/}.
The underlying algorithm is not very conceptual. It is based on
a case by case analysis of the structure of the coefficients
in the recurrence relations for the various families in the Askey scheme.
Another approach, more conceptual but not yet covering everything,
was proposed and implemented by Koepf \& Schmersau \cite{20}.
In view of the present paper, the approach of {\tt rec2ortho}, which has
to deal with 13 families in order to cover the whole Askey scheme, might be
adapted by looking for a match with one of the five cases in sections
\ref{48} and \ref{49}.

The contents of this paper are as follows. Sections 2 and 3
illustrate the general principle for the case of the classical orthogonal
polynomials. Manifolds with corners are briefly introduced in section 4.
The Askey scheme is presented in section 5. Next a summary of the results,
including the various parts of the Askey scheme covered by the five
local charts, is given in section 6. The detailed results are in sections
7 and 8, for the Racah and Wilson manifolds, respectively. In conclusion,
section 9 discusses the results and formulates work yet to be done.
\section{The general principle}
Consider the classical orthogonal polynomials (see
Chihara \cite[Ch.~V, \S2]{6}) as monic polynomials
$p_n(x)=x^n+$ terms of degree less than $n$:
\begin{itemize}
\item
Jacobi polynomials $p_n^{(\al,\be)}(x)$ with weight function
$(1-x)^\al\,(1+x)^\be$ on $(-1,1)$;
\item
Laguerre polynomials $\ell_n^\al(x)$ with weight function $e^{-x}\,x^\al$
on $(0,\iy)$;
\item
Hermite polynomials $h_n(x)$ with weight function $e^{-x^2}$ on $(-\iy,\iy)$.
\end{itemize}

They are connected by limit relations (see \cite[Ch.~2]{7}):
\begin{align}
\lim_{\be\to\iy}(-\be/2)^n\,p_n^{(\al,\be)}(1-2x/\be)&=\ell_n^\al(x),
\label{10}\\
\lim_{\al\to\iy}\al^{n/2}p_n^{(\al,\al)}(x/\al^{1/2})&=h_n(x),
\label{11}\\
\lim_{\al\to\iy} (2\al)^{-n/2}\,\ell_n^\al\big((2\al)^{1/2} x+\al\big)&=h_n(x).
\label{6}
\end{align}

The limit relation \eqref{10} is immediate from the explicit expressions
of Jacobi and Laguerre polynomials as terminating hypergeometric series
(see \cite[(1.8.1), (1.8.4), (1.11.1), (1.11.4)]{7}):
\begin{align*}
p_n^{(\al,\be)}(1-2x)&=\frac{2^n\,(\al+1)_n}{(n+\al+\be+1)_n}\,
\hyp21{-n,n+\al+\be+1}{\al+1}x,\\
\ell_n^\al(x)&=(-1)^n\,(\al+1)_n\,\hyp11{-n}{\al+1}x.
\end{align*}
Then \eqref{11} is equivalent by quadratic transformations to the cases
$\al=\pm\frac12$ of \eqref{10}.

The limit relation \eqref{6}, first observed in 1939 by Palam\`a \cite{10}
and by Toscano \cite{11},
cannot be easily obtained from explicit power series.
An approach to prove \eqref{6} by taking limits of the corresponding
weight functions was given by Askey \cite{8}. His approach works also
for \eqref{10} and \eqref{11}. 
Indeed, the limit relations for the weight functions corresponding to
\eqref{10}--\eqref{6} are:
\begin{align}
\lim_{\be\to\iy} x^\al\,(1-x/\be)^\be&=x^\al\,e^{-x},
\label{7}\\
\lim_{\al\to\iy}(1-x^2/\al)^\al&=e^{-x^2},
\label{8}\\
\lim_{\al\to\iy}\big(1+(2/\al)^{1/2}x\big)^\al\,e^{-(2\al)^{1/2} x}&=e^{-x^2}.
\label{9}
\end{align}
The limits \eqref{7} and \eqref{8} are standard limits, while \eqref{9},
observed by Askey \cite{8}, is an easy exercise.
What we in fact need in order to conclude rigorously that
\eqref{10}--\eqref{6} follow from \eqref{7}--\eqref{9}, is not just
the pointwise limits for
the corresponding weight functions, but the limits for the corresponding
moments:
\begin{proposition}
Let be given monic orthogonal polynomials $(p_n)_{n-0,1,2,\ldots}$
with respect to an orthogonality measure having moments
$\mu_n$ {\rm ($n=0,1,2,\ldots$)} and, for $\al>0$,
monic orthogonal polynomials $(p_n^\al)_{n-0,1,2,\ldots}$ with respect to
an orthogonality measure having moments
$\mu_n^\al$ {\rm ($n=0,1,2,\ldots$)}. 
If $\lim_{\al\to\iy}\mu_n^\al=\mu_n$ for all $n$
then $\lim_{\al\to\iy}p_n^\al(x)=p_n(x)$ for all $n$.
\end{proposition}
\Proof
Use that
\[
p_n(x)=\Big(\det(\mu_{i+j})_{i,j=0,1,\ldots,n-1}\Big)^{-1}
\begin{vmatrix}
\mu_0&\mu_1&\mu_2&\cdots&\mu_n\\
\mu_1&\mu_2&\mu_3&\cdots&\mu_{n+1}\\
\cdots&\cdots&\cdots&\cdots&\cdots\\
\mu_{n-1}&\mu_n&\mu_{n+1}&\cdots&\mu_{2n-1}\\
1&x&x^2&\cdots&x^n
\end{vmatrix}
\]
(see \cite[(2.2.6)]{9}), and similarly for $p_n^\al(x)$.\qed
\bLP\indent
In order to conclude from pointwise limits of weight functions as
in \eqref{7}--\eqref{9} that the corresponding limits for the moments hold,
one needs Lebesgue's dominated convergence theorem.
The relevant inequalities corresponding to \eqref{7}--\eqref{9} are:
\begin{align*}
x^\al\,(1-x/\be)^\be&\le x^\al\,e^{-x},\\
(1-x^2/\al)^\al&\le e^{-x^2},\\
\big(1+(2/\al)^{1/2}x\big)^\al\,e^{-(2\al)^{1/2} x}&\le
\begin{cases}e^{-x^2}&(-(\al/2)^{1/2}<x\le0),\\
(1+2x)\,e^{-x}&(x\ge0,\;\al\ge\thalf).
\end{cases}
\end{align*}

In this paper I want to advertize another method to prove limit formulas
of the above type for orthogonal polynomials, namely by using
the three-term recurrence relation.
The celebrated Favard theorem (see Chihara \cite[Ch.~I, Theorem 4.4]{6})
states that $\{p_n\}_{n=0,1,2,\ldots}$
is a system of monic orthogonal polynomials with respect to
a positive orthogonality
measure if and only if a recurrence relation
\begin{align}
x\,p_n(x)&=p_{n+1}(x)+B_n\,p_n(x)+C_n\,p_{n-1}(x),
\quad n=1,2,\ldots,
\label{3}\\
x\,p_0(x)&=p_1(x)+B_0\,p_0(x),
\label{4}\\
p_0(x)&=1,
\nonu
\end{align}
is valid with $C_n>0$ and $B_n$ real.
Below, when we will give the recurrence relation \eqref{3} with explicit
coefficients $B_n$ and $C_n$ depending analytically on $n$,
then we will silently assume that it also implies
the case $n=0$, i.e.\ \eqref{4}. Just take $B_n$ for $n=0$
and omit the term $C_n\,p_{n-1}(x)$ for $n=0$.

If the coefficients $B_n$ and $C_n$ are given then $p_n$ is completely
determined by this recurrence relation. In particular, if $B_n$ and $C_n$
continuously depend on some parameter $\la$ then $p_n(x)$ will also
continuously depend on $\la$.
For example, monic Hermite polynomials $h_n$
satisfy the recurrence relation
\begin{equation}
x\,h_n(x)=h_{n+1}(x)+\thalf n\,h_{n-1}(x).
\label{1}
\end{equation}
and monic Laguerre polynomials $\ell_n^\al$ satisfy the recurrence relation
\begin{equation}
x\,\ell_n^\al(x)=p_{n+1}(x)+(2n+\al+1)\,\ell_n^\al(x)+
n\,(n+\al)\,\ell_{n-1}^\al(x).
\label{5}
\end{equation}
Now consider rescaled monic Laguerre polynomials
\begin{equation*}
p_n(x)=p_n(x;\al,\rho,\si):=\rho^n\,\ell_n^\al(\rho^{-1}x-\si).
\end{equation*}
By \eqref{5} these satisfy the recurrence relation
\begin{equation}
x\,p_n(x)=p_{n+1}(x)+\rho\,(2n+\al+1+\si)\,p_n(x)+
\rho^2\,n\,(n+\al)\,p_{n-1}(x).
\label{2}
\end{equation}
We want to rescale in such a way that, as $\al\to\iy$,
$p_n(x)$ will tend to $h_n(x)$.
It is easy to see how to do this when we compare \eqref{1} and \eqref{2}.
Put $\rho:=(2\al)^{-1/2}$, $\si:=-\al$. Then \eqref{2} becomes
\begin{equation*}
x\,p_n(x)=
p_{n+1}(x)+(2\al)^{-1/2}\,(2n+1)\,p_n(x)+
\frac{n\,(n+\al)}{2\al}\,p_{n-1}(x).
\end{equation*}
The recurrence coefficients now tend to 0 resp.\ $n/2$ as $\al\to\iy$.
Hence $p_n(x)\to h_n(x)$ as $\al\to\iy$, i.e., we have recovered \eqref{6}.

\section{Uniform limit of Jacobi polynomials}
\label{21}
It is now natural to conjecture that we might also make these limit transitions
in the parameter plane in a more uniform way, i.e., to make such a rescaling
of the Jacobi polynomials that they depend continuously on $(\al,\be)$ in the
extended parameter plane and reduce to (possibly rescaled) Laguerre and
Hermite polynomials on the boundary lines and boundary vertex at infinity,
respectively.
For this purpose we consider Jacobi polynomials with arbitrary rescaling:
\begin{equation}
p_n(x):=\rho^n\,p_n^{(\al,\be)}(\rho^{-1}x-\si).
\label{12}
\end{equation}
These polynomials satisfy recurrence relations \eqref{3} with
(see \cite[(1.8.4)]{7}):
\begin{align}
C_n&:=
\rho^2\,\frac{4n\,(n+\al)\,(n+\be)\,(n+\al+\be)}
{(2n+\al+\be-1)\,(2n+\al+\be)^2\,(2n+\al+\be+1)}
\nonu\\
&\;=\frac{\rho^2\,\al\,\be}{(\al+\be)^3}\,
\frac{4n\,(1+n/\al)\,(1+n/\be)\,(1+n/(\al+\be))}
{(1+(2n-1)/(\al+\be))\,(1+2n/(\al+\be))^2\,(1+(2n+1)/(\al+\be))}
\label{13}\\
\intertext{and}
B_n&:=\rho\,\Big(\frac{\be^2-\al^2}{(2n+\al+\be)\,(2n+\al+\be+2)}+\si\Big)
\nonu\\
&\;=\rho\,\Bigl(\frac{\be-\al}{\be+\al}\;
\frac1{1+2n/(\al+\be)}\;\frac1{1+(2n+2)/(\al+\be)}+\si\Bigr).
\label{14}
\end{align}
{}From \eqref{13} we see that the choice
\begin{equation}
\rho:=\frac{(\al+\be)^{3/2}}{\al^{1/2}\,\be^{1/2}}
\label{15}
\end{equation}
makes $C_n$ continuous in $(\al,\be)$ on the extended parameter plane.
Next we see from \eqref{14} that the choice
\begin{equation}
\si:=\frac{\al-\be}{\al+\be}
\label{16}
\end{equation}
makes $B_n$ continuous in $(\al,\be)$ (extended) as well. Indeed,
we can now rewrite
\begin{align}
C_n&=\frac{4n\,(1+n/\al)\,(1+n/\be)\,(1+n/(\al+\be))}
{(1+(2n-1)/(\al+\be))\,(1+2n/(\al+\be))^2\,(1+(2n+1)/(\al+\be))}\,,
\label{19}\\
B_n&=\frac{\be^{-1}-\al^{-1}}{(\al^{-1}+\be^{-1})^{1/2}}\;
\frac{4n+2+4n\,(n+1)/(\al+\be)}
{\big(1+2n/(\al+\be)\big)\,\big(1+(2n+2)/(\al+\be)\big)}\,,
\label{20}
\end{align}
which are continuous in $(\al^{-1},\be^{-1})$ for $\al^{-1},\be^{-1}\ge0$.

As a result we can consider
the $(\al^{-1},\be^{-1})$-parameter plane.
For $\al^{-1},\be^{-1}>0$ we have the rescaled Jacobi polynomials
\eqref{12} with $\rho$ and $\si$ given by \eqref{15} and \eqref{16}.
These polynomials extend continuously to the closure
$\{(\al^{-1},\be^{-1})\mid \al^{-1},\be^{-1}\ge0\}$.
On the boundary lines $\{(\al^{-1},0)\mid\al^{-1}>0\}$
and $\{(0,\be^{-1})\mid\be^{-1}>0\}$
these polynomials, in view of \eqref{2}, become rescaled Laguerre polynomials
\begin{equation}
\rho^n\,\ell_n^\al(\rho^{-1}x-\si)\quad\rm{with}\quad
\rho=-2\al^{-1/2},\;\si=-\al,\quad\rm{resp.}\quad
\rho=2\be^{-1/2},\;\si=-\be.
\label{23}
\end{equation}
On the boundary vertex $(0,0)$ the polynomials, in view of
\eqref{1}, become rescaled Hermite polynomials
\begin{equation}
\rho^n\,h_n(\rho^{-1}x),\quad\rm{with}\quad\rho=2^{3/2}.
\label{24}
\end{equation}
\begin{remark}
Our limit
\begin{equation*}
\lim_{\al,\be\to\iy}\rho^n\,p_n^{(\al,\be)}(\rho^{-1}\,x-\si)=
2^{3n/2}\,h_n(2^{-3n/2}\,x)
\end{equation*}
with $\rho,\si$ given by \eqref{15}, \eqref{16}, implies the limit given
in \cite[\S2.6.4]{3}.
\end{remark}
\section{Manifolds with corners}
The closed subset
\begin{equation}
\{(\al^{-1},\be^{-1})\in\RR^2\mid \al^{-1},\be^{-1}\ge0\}
\label{22}
\end{equation}
of the $(\al^{-1},\be^{-1})$-parameter plane considered at the end of
\S\ref{21} is a prototype of a so-called
 manifold with corners.
More generally, define
\[
\RR_{(q)}^n:=\{(x_1,\ldots,x_n)\in\RR^n\mid x_{q+1},\ldots,x_n\ge0\}\quad
(q=0,1,\ldots,n).
\]
{\em Manifolds with corners}, introduced by Cerf \cite{12}
and Douady \cite{13}, are topological Hausdorff spaces
which are locally homeomorphic with open subsets of spaces $\RR_{(q)}^n$.
The definition is analogous to the definition of an ordinary manifold,
but there the local homeomorphisms only map onto open subsets of $\RR^n$.
Thus, for a manifold with corners denoted by $X$ we have charts
$(U,\phi)$ such that $\phi\colon U\to\phi(U)$ is a homeomorphism from
an open subset $U$ of $X$ onto an open subset $\phi(U)$ of some
$\RR_{(q)}^n$. If $(U,\phi)$ and $(V,\psi)$ are two charts on $X$ then
$\psi\circ\phi^{-1}\colon \phi(U\cap V)\to\psi(U\cap V)$ must be a
homeomorphism. If these maps extend on a larger open
subset of $\RR^n$ to $C^k$-diffeomorphisms (or $C^\iy$-diffeomorphisms
or analytic diffeomorphisms), then $X$ is called a $C^k$ (or $C^\iy$
or analytic) manifold with corners.

On the prototypical manifold with corners given by \eqref{22} each point
can be associated with a system of orthogonal polynomials. The coefficients
$C_n$ and $B_n$ given by \eqref{19} and \eqref{20} are continuous on this
manifold, and thus the resulting orthogonal polynomials $p_n$ obtained
from $C_n$ and $B_n$ by the recurrence relation \eqref{3} are continuous
on the manifold. In the interior, on the boundary lines and on the
boundary vertex these polynomials respectively become (rescaled)
Jacobi, Laguerre and Hermite polynomials as described at the end of
\S\ref{21}. Note however that $B_n$, given by \eqref{20}, is
not differentiable at the point $(\al^{-1},\be^{-1})=(0,0)$.
\section{The Askey scheme}
\label{36}
\begin{figure}[t]
\centering
\begin{picture}(300,220)
\setlength{\unitlength}{3.5mm}
\put(8.5,20) {\framebox(4.5,2) {Wilson}}
\put(9.5,20) {\vector(-1,-2){1}}
\put(12.5,20) {\vector(0,-1){2}}
\put(18.5,20) {\framebox(4.5,2) {Racah}}
\put(19.5,20) {\vector(0,-1){3}}
\put(22.5,20) {\vector(1,-2){1.5}}
\put(2.5,15) {\framebox(7,3) {\shortstack {Continuous\\[1mm]dual Hahn}}}
\put(6.5,15) {\vector(0,-1){2}}
\put(10.5,15) {\framebox(6,3) {\shortstack {Continuous\\[1mm]Hahn}}}
\put(12.5,15) {\vector(-5,-2){4.8}}
\put(13.5,15) {\vector(0,-1){3}}
\put(17.8,15) {\framebox(3.5,2) {Hahn}}
\put(18.5,15) {\vector(-4,-3){4}}
\put(19.5,15) {\vector(0,-1){3}}
\put(20.5,15) {\vector(4,-3){3.9}}
\put(22.5,15) {\framebox(6,2) {Dual Hahn}}
\put(24.5,15) {\vector(-4,-3){4}}
\put(26.5,15) {\vector(0,-1){3}}
\put(3,10) {\framebox(6,3)
{\shortstack {Meixner-\\[1mm]Pollaczek}}}
\put(7,10) {\vector(1,-2){1.5}}
\put(10.5,10) {\framebox(4.5,2) {Jacobi}}
\put(12.0,10) {\vector(0,-1){3}}
\put(14.5,10) {\vector(0,-1){8}}
\put(16.5,10) {\framebox(4.5,2) {Meixner}}
\put(17.5,10) {\vector(-5,-3){5}}
\put(19.5,10) {\vector(0,-1){3}}
\put(22.5,10) {\framebox(6.5,2) {Krawtchouk}}
\put(25,10) {\vector(-5,-3){5}}
\put(8,5) {\framebox(5,2) {Laguerre}}
\put(10.5,5) {\vector(2,-3){2}}
\put(17,5) {\framebox(4.5,2) {Charlier}}
\put(19.5,5) {\vector(-3,-2){4.5}}
\put(12.0,0) {\framebox(4.5,2) {Hermite}}
\end{picture}
\caption{Askey scheme}
\label{fig:1}
\end{figure}
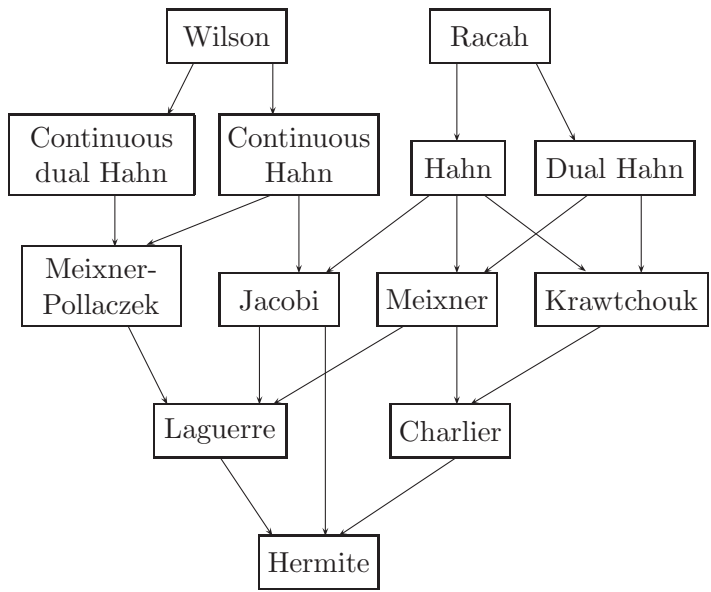

The Askey scheme is given by  Figure \ref{fig:1}.
The various families of orthogonal polynomials mentioned here are
all of classical type, i.e., the orthogonal polynomials
$\{p_n\}_{n=0,1,\ldots}$
satisfy an equation of the form
\[
L\,p_n=\la_n\,p_n,
\]
where $L$ is some second order operator (differential or difference)
which does not depend on $n$.
The arrows in the chart denote limit transitions between the various families.
The number of additional parameters on which the polynomials depend,
decreases as we go further down in the chart.
In the top row there are 4 parameters. In each subsequent row there is
one parameter less. The Hermite polynomials in the bottom row no longer
depend on parameters. The families in the left part of the chart
consist of polynomials being orthogonal \wrt an absolutely continuous
measure, while the ones displayed to the right of Hermite
are orthogonal \wrt a discrete
measure. In the case of Racah, Hahn, dual Hahn and Krawtchouk polynomials
the support
of the measure has finite cardinality, say $N+1$, and we consider
only polynomials up to degree $N$.
All the polynomials in the Askey scheme have explicit expressions as
hypergeometric functions, see \cite[Chapter 1]{7}.

{\em Racah polynomials} (see \cite[\S1.2]{7}) are defined by
\begin{equation}
R_n\big(y(y+\ga+\de+1);\al,\be,\ga,\de\big):=
\hyp43{-n,n+\al+\be+1,-y,y+\ga+\de+1}{\al+1,\be+\de+1,\ga+1}1,
\label{26}
\end{equation}
where $\ga+1=-N$ and $n=0,1,\ldots,N$.
Then we write the {\em monic Racah polynomials} as
\begin{equation}
r_n(x)=r_n(x;\al,\be,-N-1,\de):=
\frac{(\al+1)_n\,(\be+\de+1)_n\,(-N)_n}{(n+\al+\be+1)_n}\,
R_n(x;\al,\be,-N-1,\de).
\label{27}
\end{equation}

{\em Wilson polynomials} (see \cite[\S1.1]{7}) are defined by
\begin{equation}
W_n(y^2;a,b,c,d):=(a+b)_n\,(a+c)_n\,(a+d)_n\,
\hyp43{-n,n+a+b+c+d-1,a+ix,a-ix}{a+b,a+c,a+d}1.
\label{28}
\end{equation}
Then we write the {\em monic Wilson polynomials} as
\begin{equation}
w_n(x)=w_n(x;a,b,c,d):=
\frac{(-1)^n}{(n+a+b+c+d-1)_n}\,W_n(x;a,b,c,d).
\label{29}
\end{equation}

The following list
gives the notation to be used in this paper
for the monic orthogonal polynomials in the other
families (apart from Racah and Wilson) in the Askey scheme.
In each case the monic polynomials are expressed in terms of the
polynomials in usual notation and normalization as given under the heading
``Definition'' in each subsection of \cite[Chapter 1]{7}. The
constants occurring in the definitions of the monic polynomials are
taken from
the formulas for $p_n$ under the heading ``Normalized recurrence relation''
in each subsection of \cite[Chapter 1]{7}.
\mLP
{\bf Hahn} \cite[\S1.5]{7}:
\begin{equation*}
q_n(x;\al,\be,N):=\frac{(\al+1)_n\,(-N)_n}{(n+\al+\be+1)_n}\,
Q_n(x;\al,\be,N).
\end{equation*}
{\bf Dual Hahn} \cite[\S1.6]{7}:
\begin{equation*}
r_n^{\rm DH}(x;\ga,\de,N):=(\ga+1)_n\,(-N)_n\,R_n(x;\ga,\de,N).
\end{equation*}
{\bf Jacobi} \cite[\S1.8]{7}:
\begin{equation*}
p_n^{(\al,\be)}(x):=\frac{2^n\,n!}{(n+\al+\be+1)_n}\,P_n^{(\al,\be)}(x).
\end{equation*}
{\bf Meixner} \cite[\S1.9]{7}:
\begin{equation*}
m_n(x;\be,c):=(\be)_n\,\Big(\frac c{c-1}\Big)^n\,M_n(x;\be,c).
\end{equation*}
{\bf Krawtchouk} \cite[\S1.10]{7}:
\begin{equation*}
k_n(x;p,N):=(-N)_n\,p^n\,K_n(x;p,N).
\end{equation*}
{\bf Laguerre} \cite[\S1.11]{7}:
\begin{equation*}
\ell_n^\al(x):=(-1)^n\,n!\,L_n^{(\al)}(x).
\end{equation*}
{\bf Charlier} \cite[\S1.12]{7}:
\begin{equation*}
c_n(x;a):=(-1)^n\,a^n\,C_n(x;a).
\end{equation*}
{\bf Hermite} \cite[\S1.13]{7}:
\begin{equation*}
h_n(x):=2^{-n}\,H_n(x).
\end{equation*}
{\bf Continuous Hahn} \cite[\S1.4]{7}:
\begin{equation*}
p_n^{\rm CH}(x;a,b,c,d):=\frac{n!}{(n+a+b+c+d-1)_n}\,p_n(x;a,b,c,d).
\end{equation*}
{\bf Continuous dual Hahn} \cite[\S1.3]{7}:
\begin{equation*}
s_n(x;a,b,c):=(-1)^n\,S_n(x;a,b,c).
\end{equation*}
{\bf Meixner-Pollaczek} \cite[\S1.7]{7}
\begin{equation*}
p_n^{(\la)}(x;\phi):=\frac{n!}{(2\sin\phi)^n}\,P_n^{(\la)}(x;\phi).
\end{equation*}
\section{Uniform limits in the Askey scheme}
\label{25}
The monic Racah polynomials \eqref{27}
satisfy (see \cite[(1.2.4)]{7}) the three-term recurrence relation
\begin{equation*}
x\,r_n(x)=r_{n+1}(x)+B_n\,r_n(x)+C_n\,r_{n-1}(x),
\end{equation*}
where $B_n=a_n+c_n$, $C_n=a_{n-1}\,c_n$ with
\begin{align*}
a_n&:=\frac{(n+\al+1)(n+\al+\be+1)(n+\be+\de+1)(N-n)}
{(2n+\al+\be+1)(2n+\al+\be+2)}\,,\sLP
c_n&:=\frac{n(n+\al+\be+N+1)(\de-\al-n)(n+\be)}
{(2n+\al+\be)(2n+\al+\be+1)}\,.
\end{align*}
The Racah polynomials $\{r_n\}_{n=0,1,\ldots N}$ will be orthogonal with
respect to certain positive weights on some set of $N+1$ points in $\RR$ iff
$C_n>0$ for $n=1,2,\ldots,N$. A sufficient condition for this is that
the parameters satisfy the inequalities
$\al,\be>-1$, $\de>\al+N$ and $N>0$ ($N$ integer).
It will be convenient to restrict to the smaller parameter region
\begin{equation}
\al,\be>0,\quad
N>1,\quad
\de>\al+N,
\label{18}
\end{equation}
and to drop the assumption that $N$ is integer. All limits to be considered
will be in this parameter region.

Each point in the four-dimensional $(\al,\be,\de,N)$-space satisfying
\eqref{18}
thus corresponds with a system of Racah polynomials being orthogonal
with respect to a positive measure (if $N$ is moreover integer).
This defines an open subset $X_0$ of $\RR^4$.
Analogous to what we did in \S\ref{21} for an open part of
$(\al,\be)$-parameter space corresponding with Jacobi polynomials,
we want to find suitable local coordinates (charts $(U,\phi)$)
in the four-dimensional
Racah parameter space $X_0$ and suitable rescalings of the Racah
polynomials which depend on these local coordinates such that the
resulting coefficients $B_n$ and $C_n$ in the recurrence relation
have a continuous extension to the boundary of $\phi(U)$, and such
that their restriction to (parts of certain dimension of) this
boundary corresponds to certain arrows lying somewhere under the
Racah box in Figure \ref{fig:1}. $N^{-1}$ should be a monomial in these
coordinates. Thus, on approaching the boundary,  $N^{-1}$ will be either
fixed or it will tend to 0. As a consequence, it is only a minor
shortcoming in our set-up that $N$ is discrete.

\begin{figure}[t]
\hskip-2cm
\setlength{\unitlength}{3mm}
\begin{minipage}{6cm}
\begin{picture}(0,22.5)
\put(16.5,20) {\framebox(5,2) {Racah}}
\put(18.5,20) {\vector(0,-1){3}}
\put(16.5,15) {\framebox(4,2) {Hahn}}
\put(17.5,15) {\vector(-1,-1){2.5}}
\put(18.5,15) {\vector(0,-1){2.5}}
\put(19.5,15) {\vector(3,-2){3.75}}
\put(10.5,10) {\framebox(5,2.5) {Jacobi}}
\put(12.0,10) {{\vector(0,-1){2.5}}}
\put(14.5,10) {{\vector(0,-1){7.5}}}
\put(16.5,10) {\framebox(5.5,2.6) {Meixner}}
\put(17.5,10) {\vector(-5,-3){4}}
\put(19.5,10) {\vector(0,-1){2.5}}
\put(22.5,10) {\framebox(8,2.5) {Krawtchouk}}
\put(25,10) {\vector(-1,-2){1.25}}
\put(8,5) {\framebox(6,2.5) {Laguerre}}
\put(10.5,5) {\vector(1,-1){2.5}}
\put(18.5,5) {\framebox(5.5,2.5) {Charlier}}
\put(21.5,5) {\vector(-2,-1){5}}
\put(12.0,0) {\framebox(5.5,2.5) {Hermite}}
\end{picture}
\end{minipage}
\hskip-1cm
\begin{minipage}{4cm}
\begin{picture}(0,22.5)
\put(18.5,20) {\framebox(5,2.5) {Racah}}
\put(21.5,20) {\vector(0,-1){2.5}}
\put(18.5,15) {\framebox(7,2.5) {Dual Hahn}}
\put(21.5,15) {\vector(-1,-1){2.5}}
\put(21.5,15) {\vector(1,-1){2.5}}
\put(16.5,10) {\framebox(5.5,2.5) {Meixner}}
\put(19.5,10) {\vector(0,-1){2.5}}
\put(22.5,10) {\framebox(8,2.5) {Krawtchouk}}
\put(25,10) {\vector(-3,-4){1.9}}
\put(18.5,5) {\framebox(5.5,2.5) {Charlier}}
\put(21.5,5) {\vector(0,-1){2.5}}
\put(18.5,0) {\framebox(5.5,2.5) {Hermite}}
\end{picture}
\end{minipage}
\hskip0.5cm
\begin{minipage}{2cm}
\begin{picture}(0,22.5)
\put(18.5,20) {\framebox(5,2.5) {Racah}}
\put(21.5,20) {\vector(0,-1){2.5}}
\put(17.5,15) {\framebox(7,2.5) {Dual Hahn}}
\put(21.5,15) {\vector(0,-1){2.5}}
\put(18.5,10) {\framebox(5.5,2.5) {Meixner}}
\put(21.5,10) {\vector(0,-1){2.5}}
\put(18.5,5) {\framebox(6,2.5) {Laguerre}}
\put(21.5,5) {\vector(0,-1){2.5}}
\put(18.5,0) {\framebox(5.5,2.5) {Hermite}}
\end{picture}
\end{minipage}
\caption{Parts of the Racah manifold covered by the three
specific charts}
\label{fig:2}
\end{figure}
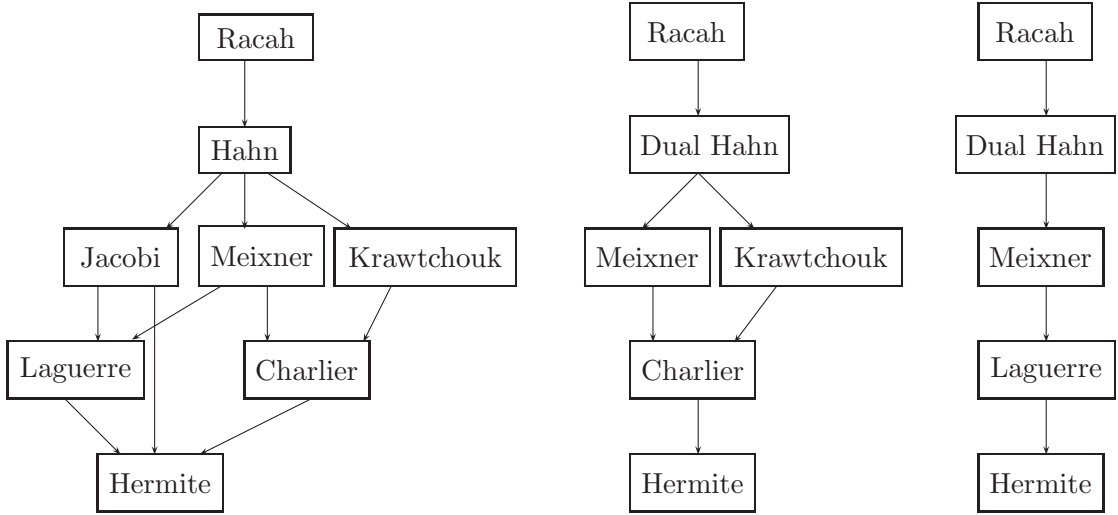

It will turn out that in this way we can extend the Racah parameter space
$X_0$ to a manifold with corners $X$ (the {\em Racah manifold})
covering everything which lies below
the Racah box in the Askey scheme if we use three local charts.
The three parts of the Askey scheme which are respectively covered
by these charts are given in Figure \ref{fig:2}.

\begin{figure}[t]
\hskip1cm
\setlength{\unitlength}{3mm}
\begin{minipage}{5cm}
\begin{picture}(0,22.5)
\put(10.5,20) {\framebox(5,2.5) {Wilson}}
\put(13.1,20) {\vector(0,-1){1.5}}
\put(10.5,15) {\framebox(5,3.5) {\shortstack {Cont.\\[1mm]Hahn}}}
\put(12.5,15) {\vector(-4,-3){3}}
\put(13.5,15) {\vector(0,-1){2.5}}
\put(3,10) {\framebox(6.5,3.5)
{\shortstack {Meixner-\\[1mm]Pollaczek}}}
\put(7.5,10) {\vector(1,-2){1.3}}
\put(10.5,10) {\framebox(5,2.5) {Jacobi}}
\put(12.0,10) {\vector(0,-1){2.5}}
\put(14.5,10) {\vector(0,-1){7.5}}
\put(8,5) {\framebox(6,2.5) {Laguerre}}
\put(11.5,5) {\vector(1,-2){1.2}}
\put(12.0,0) {\framebox(5.5,2.5) {Hermite}}
\end{picture}
\end{minipage}
\hskip1.5cm
\begin{minipage}{3cm}
\begin{picture}(0,22.5)
\put(1,20) {\framebox(5,2.5) {Wilson}}
\put(3.5,20) {\vector(0,-1){1.5}}
\put(0,15) {\framebox(7,3.5) {\shortstack {Cont.\\[1mm]dual Hahn}}}
\put(3.5,15) {\vector(0,-1){1.5}}
\put(0.25,10) {\framebox(6.5,3.5)
{\shortstack {Meixner-\\[1mm]Pollaczek}}}
\put(3.5,10) {\vector(0,-1){2.5}}
\put(0.5,5) {\framebox(6,2.5) {Laguerre}}
\put(3.5,5) {\vector(0,-1){2.5}}
\put(0.75,0) {\framebox(5.5,2.5) {Hermite}}
\end{picture}
\end{minipage}
\hskip1cm
\begin{minipage}{3cm}
\begin{picture}(0,22.5)
\put(1,20) {\framebox(5,2.5) {Wilson}}
\put(3.5,20) {\vector(0,-1){1.5}}
\put(0,15) {\framebox(7,3.5) {\shortstack {Cont.\\[1mm]dual Hahn}}}
\put(3.5,15) {\vector(0,-1){7.5}}
\put(0.5,5) {\framebox(6,2.5) {Laguerre}}
\put(3.5,5) {\vector(0,-1){2.5}}
\put(0.75,0) {\framebox(5.5,2.5) {Hermite}}
\end{picture}
\end{minipage}
\caption{Limits of Wilson polynomials: on the left parameters in two
complex conjugate pairs; in the middle one complex conjugate pair and
two real; on the right all parameters real}
\label{fig:3}
\end{figure}
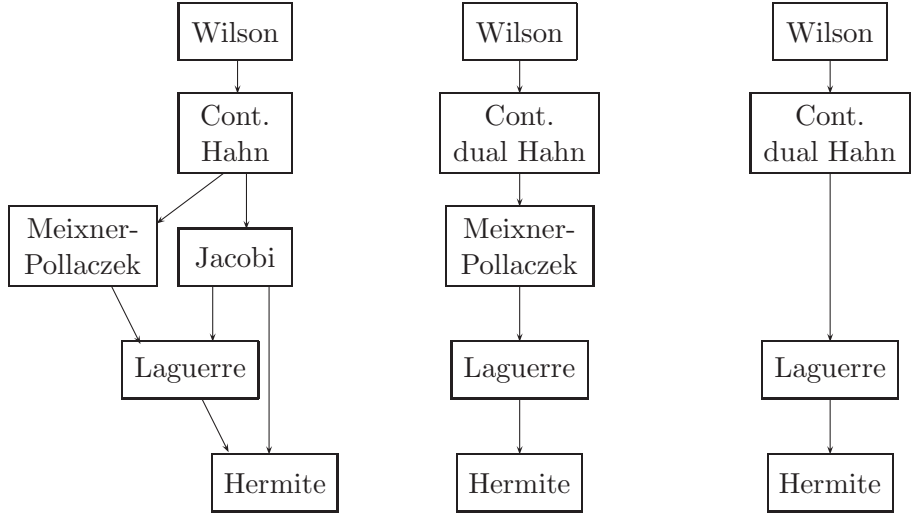

The monic Wilson polynomials \eqref{29} satisfy
(see \cite[(1.1.5)]{7}) the three-term recurrence relation
\begin{equation*}
x\,w_n(x)=w_{n+1}(x)+B_n\,w_n(x)+C_n\,w_{n-1}(x),
\end{equation*}
where $B_n=a_n+c_n-a^2$, $C_n=a_{n-1}\,c_n$ with
\begin{align*}
a_n&:=\frac{(n+a+b+c+d-1)\,(n+a+b)\,(n+a+c)\,(n+a+d)}
{(2n+a+b+c+d-1)\,(2n+a+b+c+d)}\,,\\
c_n&:=\frac{n\,(n+b+c-1)\,(n+b+d-1)\,(n+c+d-1)}
{(2n+a+b+c+d-2)\,(2n+a+b+c+d-1)}\,.
\end{align*}
Then $B_n$ and $C_n$ turn out to be symmetric in $a,b,c,d$, and therefore
$w_n(x;a,b,c,d)$ is symmetric in $a,b,c,d$.
The polynomials $\{w_n\}_{n=0,1,2,\ldots}$ will be orthogonal with
respect to a certain positive measure on $\RR$ iff
$C_n>0$ for $n=1,2,\ldots\;$. This, in its turn, is equivalent to one
of the following three cases:
\begin{enumerate}
\item
$a,b,c,d$ are non-real and occur in complex conjugate pairs with
positive real parts.
\item
Two of the four parameters $a,b,c,d$ are non-real and form a complex
conjugate pair with positive real part.
The other two parameters are real and their sum is positive.
\item
$a,b,c,d$ are real and $C_n>0$ for $n=1,2,\ldots\;$.
\end{enumerate}
Limits of Wilson polynomials in the Askey scheme occur in all these three
cases, see Figure~\ref{fig:3}.
In case~1 we have a limit to the continuous
Hahn polynomials and to everything below the continuous Hahn polynomials.
In case 2 we have a limit to the continuous dual
Hahn polynomials and to everything below the continuous dual Hahn polynomials.
To each of these two cases will correspond one manifold with corners.
In case 3 we also have a limit to the continuous dual
Hahn polynomials, but from there a limit to the two-parameter level seems
to be missing; one rather goes to Laguerre as follows:
\begin{equation}
\lim_{c\to\iy} c^{-n}\,s_n(c\,x;a,b,c)=\ell_n^{\,a+b-1}(x).
\label{53}
\end{equation}
This is straightforward from
\cite[(1.3.1), (1.11.1)]{7} and the formulas for the monic
continuous dual Hahn and Laguerre polynomials in \S\ref{36}.
The limit \eqref{53} is missing in \cite[Appendix]{14} and
\cite[\S2.3]{7}.

It will turn out that one chart is sufficient in cases 1 and 2
in order to cover everything which is below the Wilson box. I did not
try to find a chart for case 3.

The following simple observation will be continuously used in the next
sections. If
\begin{equation}
x\,p_n(x)=p_{n+1}(x)+B_n\,p_n(x)+C_n\,p_{n-1}(x),
\label{32}
\end{equation}
and if
\begin{equation}
q_n(x):=\rho^n\,p_n(\rho^{-1}x-\si)
\label{33}
\end{equation}
then
\begin{equation*}
x\,q_n(x)=q_{n+1}(x)+\rho\,(B_n+\si)\,q_n(x)+\rho^2\,C_n\,q_{n-1}(x),
\end{equation*}

\section{The Racah manifold}
\label{48}
In this section I will present in detail the three charts covering the Racah
manifold, as introduced in \S\ref{25}, and successively corresponding to
the three graphs in Figure \ref{fig:2}.
\subsection{From Racah to Hermite along Hahn}
\label{35}
Here we will see the chart for the Racah manifold corresponding to the
first graph in Figure \ref{fig:2}.
The monic Racah polynomial $r_n$ is given by \eqref{27}.
For $t_1,t_2,t_3,t_4>0$ and $t_1\,t_3<1$, $t_2\,t_4<1$ put
\begin{align}
&p_n(x)=p_n(x;t_1,t_2,t_3,t_4):=
\rho^n\,r_n(\rho^{-1}x-\si;\al,\be,-N-1,\de),
\label{34}\sLP
{\rm where}\quad
&\al=\frac1{t_1}\,,\quad
\be=\frac1{t_1\,t_2}\,,\quad
N=\frac1{t_2\,t_4}\,,\quad
\de=\frac{1+t_2\,t_3\,t_4}{t_1\,t_2\,t_3\,t_4}\,,
\label{39}\sLP
&\rho=\frac{t_1\,t_2\,(1+t_2)^{3/2}\,t_3\,t_4^2}
{(t_1+t_4+t_2\,t_4)^{1/2}\,\big(1+(1+t_2)\,t_3\,t_4\big)^{1/2}}\,,\quad
\si=-\,\frac{(1+t_1)\,\big(1+(1+t_2+t_1\,t_2)\,t_3\,t_4\big)}
{t_1\,t_2\,(1+t_2+2\,t_1\,t_2)\,t_3\,t_4^2}\,.
\nonu
\end{align}
Then the inequalities for $t_1,t_2,t_3,t_4$ are equivalent with \eqref{18}
and the recurrence relation \eqref{32} holds with
\begin{align}
&B_n=-\,\frac{n(1+t_2)^{3/2}\,\big(1+t_2+(n+1)\,t_1\,t_2\big)}
{(1+t_2+2\,t_1\,t_2)\,(1+t_2+2\,n\,t_1\,t_2)\,
(1+t_2+2\,(n+1)\,t_1\,t_2)\,(1+t_3\,t_4+t_2\,t_3\,t_4)^{1/2}}
\nonu\\
&\times
\Big(2\,n\,t_1\,t_2\,t_3\,t_4^2\,(1+t_2+(n+1)\,t_1\,t_2)\,(1+t_2+2\,t_1\,t_2)+
t_2\,t_3\,t_4^2\,(1+t_1)\,(1+t_2)\,(1+t_2+2\,t_1\,t_2)
\nonu\\
&\qquad-t_4\,(1-t_2^2)\,(1+t_1\,t_3)
-2t_1\,(1-t_2)\,(1+t_2\,t_4)\Big)\Big/\big(t_1+t_4+t_2\,t_4\big)^{1/2},
\label{30}
\sLP
&C_n=
\frac{(1+t_2+n\,t_1\,t_2)\,\big(1+(1-n)\,t_2\,t_4)\,
(1-n\,t_1\,t_2\,t_3\,t_4)\,(1+t_3\,t_4+t_2\,t_3\,t_4+n\,t_1\,t_2\,t_3\,t_4)}
{\big(1+t_2+(2n-1)\,t_1\,t_2\big)\,(1+t_2+2\,n\,t_1\,t_2)^2\,
\big(1+t_2+(2n+1)\,t_1\,t_2\big)}
\nonu\\
&\qquad\quad\times n\,(1+n\,t_1)\,(1+n\,t_1\,t_2)\,
\frac{(1+t_2)^3}{1+t_3\,t_4+t_2\,t_3\,t_4}\,
\Big(1+(n+1)\,\frac{t_1\,t_2\,t_4}{t_1+t_4+t_2t_4}\Big).
\label{31}
\end{align}
Note that $B_n$ and $C_n$, as functions of $t_1,t_2,t_3,t_4>0$,
can be uniquely extended to continuous (but not differentiable)
functions of
$t_1,t_2,t_3,t_4\ge0$. Indeed, for $B_n$ observe that,
if $t_1,t_2,t_3,t_4>0$, then
\begin{equation*}
0\le\frac{t_1}{(t_1+t_4+t_2t_4)^{1/2}}\le t_1^{1/2},\qquad
0\le\frac{t_4}{(t_1+t_4+t_2t_4)^{1/2}}\le t_4^{1/2}.
\end{equation*}
For $C_n$ observe that,
if $t_1,t_2,t_3,t_4>0$, then
\begin{equation*}
0\le\frac{t_1t_2t_4}{t_1+t_4+t_2t_4}\le
\frac{\tfrac14 (t_1+t_2t_4)^2}{t_1+t_2t_4}=\tfrac14(t_1+t_2t_4).
\end{equation*}

We now put one or more of the $t_1,t_2,t_3,t_4$ equal to zero
in \eqref{30} and \eqref{31}, and then compare the resulting
recurrence relation \eqref{32} with the various normalized recurrence
relations in \cite[Chapter~1]{7}. In each case we find a match with
one of these normalized recurrence relations, after rescaling as in
\eqref{33} for some special $\rho$, $\si$. Then the explicitly obtained
family in the Askey scheme satisfying this recurrence relation will
be equal to the corresponding parameter restriction of the polynomial
in \eqref{34}. We obtain the following results, using the notation for
monic polynomials in the Askey scheme as given at the end of \S\ref{36}.
\mLP
{\bf Hahn}:
\begin{align}
&p_n(x;t_1,t_2,0,t_4)=
\rho^n\,q_n(\rho^{-1}x-\si;\al,\be,N),
\nonu\sLP
&\al=\frac1{t_1}\,,\quad
\be=\frac1{t_1\,t_2}\,,\quad
N=\frac1{t_2\,t_4}\,,\quad
\rho=\frac{(1+t_2)^{3/2}\,t_4}
{(t_1+t_4+t_2\,t_4)^{1/2}}\,,\quad
\si=-\,\frac{1+t_1}
{(1+t_2+2\,t_1\,t_2)\,t_4}\,.
\label{50}
\end{align}
{\bf Jacobi}:
\begin{align*}
&p_n(x;t_1,t_2,t_3,0)=p_n(x;t_1,t_2,0,0)=
\rho^n\,p_n^{(\al,\be)}(\rho^{-1}x-\si),\sLP
&\al=\frac1{t_1}\,,\quad
\be=\frac1{t_1\,t_2}\,,\quad
\rho=-\,\frac{(1+t_2)^{3/2}}
{2\,t_1^{1/2}\,t_2}\,,\quad
\si=\,\frac{-1+t_2}
{1+t_2+2\,t_1\,t_2}\,.
\end{align*}
{\bf Meixner}:
\begin{align}
&p_n(x;t_1,0,t_3,t_4)
=p_n\left(x;t_1,0,0,\frac{t_4(1-t_1\,t_3)}{1+t_3\,t_4}\right)=
\rho^n\,m_n(\rho^{-1}x-\si;\be,c),
\nonu\sLP
&\be=\frac{1+t_1}{t_1}\,,\quad
c=\frac{t_1\,(1+t_3\,t_4)}{t_1+t_4}\,,\quad
\rho=\frac{(1-t_1\,t_3)\,t_4}
{(t_1+t_4)^{1/2}\,(1+t_3\,t_4)^{1/2}}\,,\quad
\si=-\frac{(1+t_1)\,(1+t_3\,t_4)}{(1-t_1\,t_3)\,t_4}\,.
\label{52}
\end{align}
{\bf Krawtchouk}:
\begin{align*}
&p_n(x;0,t_2,t_3,t_4)=
p_n\left(x;0,t_2\,(1+t_3\,t_4+t_2\,t_3\,t_4),0,
\frac{t_4}{1+t_3\,t_4+t_2\,t_3\,t_4}\right):=
\rho^n\,k_n(\rho^{-1}x-\si;p,N),\sLP
&p=\frac{t_2\,(1+t_3\,t_4+t_2\,t_3\,t_4)}{(1+t_2)\,(1+t_2\,t_3\,t_4)}\,,\quad
N=\frac1{t_2\,t_4}\,,\\
&\rho=\frac{t_4^{1/2}\,(1+t_2)\,(1+t_2\,t_3\,t_4)}
{(1+t_3\,t_4+t_2\,t_3\,t_4)^{1/2}},\quad
\si=-\,\frac{1+t_3\,t_4+t_2\,t_3\,t_4}{t_4\,(1+t_2)\,(1+t_2\,t_3\,t_4)}\,.
\end{align*}
{\bf Laguerre}:
\begin{align*}
&p_n(x;t_1,0,t_3,0)=p_n(x;t_1,0,0,0)=
\rho^n\,\ell_n^{(\al)}(\rho^{-1}x-\si),\sLP
&\al=\frac1{t_1}\,,\quad
\rho=t_1^{1/2},\quad
\si=-\,\frac{1+t_1}{t_1}\,.
\end{align*}
{\bf Charlier}:
\begin{align*}
&p_n(x;0,0,t_3,t_4)=p_n\left(x;0,0,0,\frac{t_4}{1+t_3\,t_4}\right)=
\rho^n\,c_n(\rho^{-1}x-\si;a),\sLP
&a=\frac{1+t_3\,t_4}{t_4}\,,\quad
\rho=\frac{t_4^{1/2}}{(1+t_3\,t_4)^{1/2}}\,,\quad
\si=-\,\frac{1+t_3\,t_4}{t_4}\,.
\end{align*}
{\bf Hermite}:
\begin{align}
&p_n(x;0,t_2,t_3,0)=p_n(x;0,t_2,0,0)=p_n(x;0,0,t_3,0)=p_n(x;0,0,0,0)=
\rho^n\,h_n(\rho^{-1}x-\si),
\nonu\sLP
&\rho=2^{1/2},\quad
\si=0.
\label{51}
\end{align}

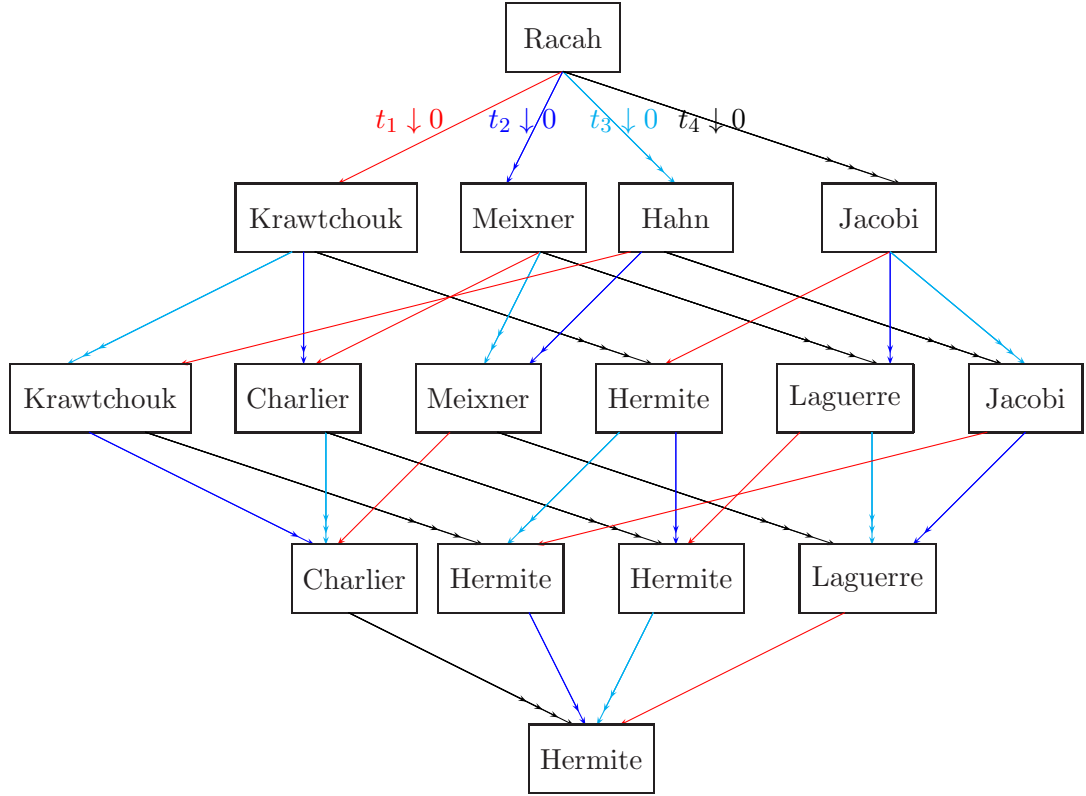
\begin{figure}[t]
\setlength{\unitlength}{3mm}
\begin{picture}(0,37.5)
\put(22,34.5) {\framebox(5,3) {Racah}}
\put(24.5,34.5) {\color{red}\vector(-2,-1){10}}
\put(16.2,32){\color{red}{$t_1\downarrow0$}}
\put(24.5,34.5) {\color{blue}\vector(-1,-2){2.5}}
\put(24.5,34.5) {\color{blue}\vector(-1,-2){2.2}}
\put(21.2,32){\color{blue}{$t_2\downarrow0$}}
\put(24.5,34.5) {\color{cyan}\vector(1,-1){5}}
\put(24.5,34.5) {\color{cyan}\vector(1,-1){4.5}}
\put(24.5,34.5) {\color{cyan}\vector(1,-1){4.0}}
\put(25.7,32){\color{cyan}{$t_3\downarrow0$}}
\put(24.5,34.5) {\color{black}\vector(3,-1){15}}
\put(24.5,34.5) {\color{black}\vector(3,-1){14}}
\put(24.5,34.5) {\color{black}\vector(3,-1){13}}
\put(24.5,34.5) {\color{black}\vector(3,-1){12}}
\put(29.6,32){\color{black}{$t_4\downarrow0$}}
\put(10,26.5) {\framebox(8,3) {Krawtchouk}}
\put(20,26.5) {\framebox(5.5,3) {Meixner}}
\put(27,26.5) {\framebox(5,3) {Hahn}}
\put(36,26.5) {\framebox(5,3) {Jacobi}}
\put(12.5,26.5) {\color{cyan}\vector(-2,-1){10}}
\put(12.5,26.5) {\color{cyan}\vector(-2,-1){9.3}}
\put(12.5,26.5) {\color{cyan}\vector(-2,-1){8.6}}
\put(13,26.5) {\color{blue}\vector(0,-1){5}}
\put(13,26.5) {\color{blue}\vector(0,-1){4.5}}
\put(13.5,26.5) {\color{black}\vector(3,-1){15}}
\put(13.5,26.5) {\color{black}\vector(3,-1){14}}
\put(13.5,26.5) {\color{black}\vector(3,-1){13}}
\put(13.5,26.5) {\color{black}\vector(3,-1){12}}
\put(23.5,26.5) {\color{red}\vector(-2,-1){10}}
\put(23.5,26.5) {\color{cyan}\vector(-1,-2){2.5}}
\put(23.5,26.5) {\color{cyan}\vector(-1,-2){2.2}}
\put(23.5,26.5) {\color{cyan}\vector(-1,-2){1.9}}
\put(23.5,26.5) {\color{black}\vector(3,-1){15}}
\put(23.5,26.5) {\color{black}\vector(3,-1){14}}
\put(23.5,26.5) {\color{black}\vector(3,-1){13}}
\put(23.5,26.5) {\color{black}\vector(3,-1){12}}
\put(27.5,26.5) {\color{red}\vector(-4,-1){20}}
\put(28,26.5) {\color{blue}\vector(-1,-1){5}}
\put(28,26.5) {\color{blue}\vector(-1,-1){4.5}}
\put(29,26.5) {\color{black}\vector(3,-1){15}}
\put(29,26.5) {\color{black}\vector(3,-1){15}}
\put(29,26.5) {\color{black}\vector(3,-1){14}}
\put(29,26.5) {\color{black}\vector(3,-1){13}}
\put(29,26.5) {\color{black}\vector(3,-1){12}}
\put(39,26.5) {\color{red}\vector(-2,-1){10}}
\put(39,26.5) {\color{blue}\vector(0,-1){5}}
\put(39,26.5) {\color{blue}\vector(0,-1){4.5}}
\put(39,26.5) {\color{cyan}\vector(6,-5){6}}
\put(39,26.5) {\color{cyan}\vector(6,-5){5.5}}
\put(39,26.5) {\color{cyan}\vector(6,-5){5}}
\put(0,18.5) {\framebox(8,3) {Krawtchouk}}
\put(10,18.5) {\framebox(5.5,3) {Charlier}}
\put(18,18.5) {\framebox(5.5,3) {Meixner}}
\put(26,18.5) {\framebox(5.5,3) {Hermite}}
\put(34,18.5) {\framebox(6,3) {Laguerre}}
\put(42.5,18.5) {\framebox(5,3) {Jacobi}}
\put(3.5,18.5) {\color{blue}{\vector(2,-1){10}}}
\put(3.5,18.5) {\color{blue}{\vector(2,-1){9.3}}}
\put(6,18.5) {\color{black}{\vector(3,-1){15}}}
\put(6,18.5) {\color{black}{\vector(3,-1){14}}}
\put(6,18.5) {\color{black}{\vector(3,-1){13}}}
\put(6,18.5) {\color{black}{\vector(3,-1){12}}}
\put(14,18.5) {\color{cyan}\vector(0,-1){5}}
\put(14,18.5) {\color{cyan}\vector(0,-1){4.5}}
\put(14,18.5) {\color{cyan}\vector(0,-1){4}}
\put(14,18.5) {\color{black}\vector(3,-1){15}}
\put(14,18.5) {\color{black}\vector(3,-1){14}}
\put(14,18.5) {\color{black}\vector(3,-1){13}}
\put(14,18.5) {\color{black}\vector(3,-1){12}}
\put(19.5,18.5) {\color{red}\vector(-1,-1){5}}
\put(21.6,18.5) {\color{black}\vector(3,-1){15}}
\put(21.6,18.5) {\color{black}\vector(3,-1){14}}
\put(21.6,18.5) {\color{black}\vector(3,-1){13}}
\put(21.6,18.5) {\color{black}\vector(3,-1){12}}
\put(27,18.5) {\color{cyan}\vector(-1,-1){5}}
\put(27,18.5) {\color{cyan}\vector(-1,-1){4.5}}
\put(27,18.5) {\color{cyan}\vector(-1,-1){4}}
\put(29.5,18.5) {\color{blue}\vector(0,-1){5}}
\put(29.5,18.5) {\color{blue}\vector(0,-1){4.5}}
\put(35,18.5) {\color{red}\vector(-1,-1){5}}
\put(38.2,18.5) {\color{cyan}\vector(0,-1){5}}
\put(38.2,18.5) {\color{cyan}\vector(0,-1){4.5}}
\put(38.2,18.5) {\color{cyan}\vector(0,-1){4}}
\put(43.3,18.5) {\color{red}\vector(-4,-1){20}}
\put(45,18.5) {\color{blue}\vector(-1,-1){5}}
\put(45,18.5) {\color{blue}\vector(-1,-1){4.5}}
\put(12.5,10.5) {\framebox(5.5,3) {Charlier}}
\put(19,10.5) {\framebox(5.5,3) {Hermite}}
\put(27,10.5) {\framebox(5.5,3) {Hermite}}
\put(35,10.5) {\framebox(6,3) {Laguerre}}
\put(15,10.5) {\color{black}\vector(2,-1){10}}
\put(15,10.5) {\color{black}\vector(2,-1){9.3}}
\put(15,10.5) {\color{black}\vector(2,-1){8.6}}
\put(15,10.5) {\color{black}\vector(2,-1){7.9}}
\put(23,10.5) {\color{blue}\vector(1,-2){2.5}}
\put(23,10.5) {\color{blue}\vector(1,-2){2.2}}
\put(28.5,10.5) {\color{cyan}\vector(-1,-2){2.5}}
\put(28.5,10.5) {\color{cyan}\vector(-1,-2){2.2}}
\put(28.5,10.5) {\color{cyan}\vector(-1,-2){1.9}}
\put(37,10.5) {\color{red}\vector(-2,-1){10}}
\put(23.0,2.5) {\framebox(5.5,3) {Hermite}}
\end{picture}
\vskip-1cm
\caption{The various parameter restrictions in the first chart for
the Racah manifold}
\label{fig:4}
\end{figure}

The various parameter restrictions of the polynomial \eqref{34}
are summarized in Figure \ref{fig:4}. Note that Racah and Hahn occur
in one box, most others in two boxes, and Hermite even in three boxes.
For instance, we reach Jacobi from Racah by putting $t_4=0$, but
then there is also no dependence on $t_3$, so we may put $t_3=0$ as well.
The Meixner and Krawtchouk cases are slightly more complicated.
For instance, we reach Meixner from Racah by putting $t_2=0$, but the
resulting expression does not change if we replace $t_3,t_4$ by
$0,t_4(1-t_1t_3)(1+t_3t_4)^{-1}$.

\begin{remark}
As a special case of the above results we have the following
limit from Hahn to Hermite:
\begin{equation*}
\lim_{t_1,t_4\downarrow0}p_n(x;t_1,t_2,0,t_4)=p_n(0,t_2,0,0).
\end{equation*}
In view of \eqref{50} and \eqref{51} we can rewrite this limit as
\begin{equation*}
\lim_{a\to\iy}\rho^n\,q_n(\rho^{-1}\,x-\si;a,a\,b,a\,N)=
2^{n/2}\,h_n(2^{-1/2}\,x),
\end{equation*}
where
\begin{equation*}
\rho=\frac{(b+1)^{3/2}}{a^{1/2}\,b^{1/2}\,N^{1/2}\,(b+N+1)^{1/2}}\,,\quad
\si=-\,\frac{a\,(a+1)\,N}{a\,b+a+2}\,.
\end{equation*}
This limit is equivalent to the limit given in \cite[(16)]{23}.

Similarly, by \eqref{52} and \eqref{51} we can rewrite the restriction
of $p_n(t_1,0,0,t_4)$ to $t_1=t_4=0$ as
\begin{equation*}
\lim_{\be\to\iy}\rho^n\,
m_n\left(\rho^{-1}\,x-\si;\frac1{\be-1}\,,\frac{1-c}{c\,(\be-1)}\right)
=2^{n/2}\,h_n(2^{-1/2}\,x),
\end{equation*}
where
\begin{equation*}
\rho=\frac{1-c}{c^{1/2}\,(\be-1)^{1/2}}\,,\quad
\si=-\,\frac{\be\,c}{1-c}\,.
\end{equation*}
This limit is equivalent to the limit given in \cite[\S2.4.4]{3}.
\end{remark}
\subsection{From Racah to Hermite along dual Hahn and Charlier}
Next we will see the chart for the Racah manifold corresponding to the
second graph in Figure \ref{fig:2}.
Let $r_n$ again be the monic Racah polynomial given by \eqref{27}.
For $s_1,s_2,s_3,s_4>0$ and $s_2^2\,s_4<1$ put
\begin{align}
&p_n(x)=p_n(x;s_1,s_2,s_3,s_4):=
\rho^n\,r_n(\rho^{-1}x-\si;\al,\be,-N-1,\de),
\label{37}\sLP
{\rm where}\quad
&\al=\frac{1+s_1}{s_1\,s_2}\,,\quad
\be=\frac{1+s_1}{s_1\,s_2^2\,s_3\,s_4}\,,\quad
N=\frac1{s_2^2\,s_4}\,,\quad
\de=\frac{1+s_1+s_2\,s_4\,(1+s_1+s_1\,s_2)}{s_1\,s_2^2\,s_4}\,,
\label{40}\sLP
&\rho=\frac{s_1\,s_2^{5/2}\,s_4}{2^{1/2}\,(1+s_1)}\,,\quad
\si=-\,\frac{(1+s_1)(1+s_3-s_2^2\,s_4)+s_1\,s_2}{s_1\,s_2^3\,s_4}\,.
\nonu
\end{align}
Then the inequalities for $s_1,s_2,s_3,s_4$ are equivalent with \eqref{18}.
Furthermore \eqref{32} holds with
\begin{align*}
&B_n=-2^{-{1/2}}\,s_2^{1/2}\,
\Big(2\,n^2\,(n+1)^2\,s_1^3\,s_2^6\,s_3^2\,s_4^3\,(1+s_1)^{-1}
+4\,n^2\,(n+1)\,s_1^2\,s_2^4\,s_3\,s_4^2\,(1+s_2\,s_3\,s_4)\\
&\qquad+n^2\,s_1\,s_2^2\,s_4\,
\big(2+2s_1-s_3-2\,s_1\,s_3-2\,s_1\,s_3^2+5\,s_2\,s_3\,s_4
+5\,s_1\,s_2\,s_3\,s_4\displaybreak[0]\\
&\qquad+s_2\,s_3^2\,s_4+2\,s_1\,s_2\,s_3^2\,s_4
+4\,s_1\,s_2\,s_3^3\,s_4
+3\,s_2^2\,s_3^2\,s_4^2+3\,s_1\,s_2^2\,s_3^2\,s_4^2
-2\,s_1\,s_2^3\,s_3^2\,s_4^2\big)
\displaybreak[0]\\
&\qquad-n\,(1+s_1+s_2\,s_3\,s_4+s_1\,s_2\,s_3\,s_4+s_1\,s_2^2\,s_3\,s_4)\,
\big(1+2\,s_1+2\,s_1\,s_3-s_2\,s_4-s_1\,s_2\,s_4\\
&\qquad\qquad
-s_2\,s_3\,s_4-2\,s_1\,s_2\,s_3\,s_4-4\,s_1\,s_2\,s_3^2\,s_4
-s_2^2\,s_3\,s_4^2
-s_1\,s_2^2\,s_3\,s_4^2+2\,s_1\,s_2^2\,s_3\,s_4^2\big)\\
&\qquad+(1+s_1)\,(1+s_2\,s_3\,s_4)\,
\big(-s_1\,s_3-s_2\,s_4-s_1\,s_2\,s_4+s_3^2\,s_4+s_1\,s_3^2\,s_4
+2\,s_1\,s_2\,s_3^2\,s_4\\
&\qquad\qquad
-s_2^2\,s_3\,s_4^2-s_1\,s_2^2\,s_3\,s_4^2
-2\,s_1\,s_2^3\,s_3\,s_4^2\big)\Big)
\displaybreak[0]\\
&\qquad\times\big((1+s_1)\,
(1+s_2\,s_3\,s_4)+2\,n\,s_1\,s_2^2\,s_3\,s_4\big)^{-1}\,
\big((1+s_1)\,(1+s_2\,s_3\,s_4)+2\,(n+1)\,s_1\,s_2^2\,s_3\,s_4\big)^{-1},
\displaybreak[0]\mLP
&C_n=n\,(1+s_1+n\,s_1\,s_2)\,\big(1+(1-n)\,s_2^2\,s_4\big)\,
\big(1+s_1+(1-n)\,s_1\,s_2^2\,s_4\big)\sLP
&\qquad\times\frac{\big(1+s_1+n\,s_1\,s_2^2\,s_3\,s_4\big)\,
\big((1+s_1)\,(1+s_2\,s_3\,s_4)+s_1\,s_3+(n+1)\,s_1\,s_2^2\,s_3\,s_4\big)}
{2\,(1+s_1)^2\,
\big((1+s_1)\,(1+s_2\,s_3\,s_4)+(2n-1)\,s_1\,s_2^2\,s_3\,s_4\big)}\sLP
&\qquad\times\frac{
\big((1+s_1)\,(1+s_2\,s_3\,s_4)+n\,s_1\,s_2^2\,s_3\,s_4\big)\,
\big((1+s_1)\,(1+s_3+s_2\,s_3\,s_4)+(n+1)\,s_1\,s_2^2\,s_3\,s_4\big)}
{\big((1+s_1)\,(1+s_2\,s_3\,s_4)+2n\,s_1\,s_2^2\,s_3\,s_4\big)^2\,
\big((1+s_1)\,(1+s_2\,s_3\,s_4)+(2n+1)\,s_1\,s_2^2\,s_3\,s_4\big)}\,.
\end{align*}
Note that $B_n$ and $C_n$, as functions of $s_1,s_2,s_3,s_4>0$,
can be uniquely extended to continuous functions of
$s_1,s_2,s_3,s_4\ge0$.

We now put one or more of the $s_1,s_2,s_3,s_4$ equal to zero and
we proceed as in \S\ref{35}. We obtain:
\mLP
{\bf dual Hahn}:
\begin{align*}
&p_n(x;s_1,s_2,0,s_4)=
\rho^n\,r_n^{\rm DH}(\rho^{-1}x-\si;\ga,\de,N),\sLP
&\ga=\frac{1+s_1}{s_1\,s_2}\,,\quad
\de=\frac1{s_1\,s_2^2\,s_4}\,,\quad
N=\frac1{s_2^2\,s_4}\,,\\
&\rho=\frac{s_1\,s_2^{5/2}\,s_4}
{2^{1/2}\,(1+s_1)}\,,\quad
\si=-\,\frac{(1+s_1)\,(1-s_2^2\,s_4)+s_1\,s_2}
{s_1\,s_2^2\,s_4}\,.
\end{align*}
{\bf Meixner}:
\begin{align*}
&p_n(x;s_1,s_2,s_3,0)
=\rho_0^n\,p_n\left(\rho_0^{-1}\,x-\si_0;s_1\,(1+s_3),
\frac{s_2(1+s_1+s_1\,s_3)}{(1+s_1)\,(1+s_3)},0,0\right)\\
&\hskip2.8cm
=\rho^n\,m_n(\rho^{-1}x-\si;\be,c),\sLP
&\be=\frac{1+s_1+s_1\,s_2}{s_1\,s_2}\,,\;\;
c=\frac{s_1\,(1+s_3)}{1+s_1+s_1\,s_3},\;\;
\rho=\frac{s_2^{1/2}}{2^{1/2}\,(1+s_1)}\,,\;\;
\si=-\,\frac{1+s_1+s_3+s_1\,s_2+s_1\,s_3}{s_2}\,,\\
&\rho_0=\frac{(1+s_3)^{1/2}\,(1+s_1+s_1\,s_3)^{1/2}}{(1+s_1)^{1/2}}\,,\quad
\si_0=\frac{s_1\,s_2^{1/2}\,s_3}
{2^{1/2}\,(1+s_1)^{1/2}\,(1+s_3)^{1/2}\,(1+s_1+s_1\,s_3)^{1/2}}\,.
\end{align*}
{\bf Krawtchouk}:
\begin{align*}
&p_n(x;0,s_2,s_3,s_4)
=\rho_0^n\,p_n\left(\rho_0^{-1}x-\si_0;
\frac{s_2}{1+s_3+s_2\,s_3\,s_4},0,s_4\,(1+s_3+s_2\,s_3\,s_4)^2\right)\\
&\hskip2.8cm=\rho^n\,k_n(\rho^{-1}x-\si;p,N),
\displaybreak[0]\sLP
&p=\frac{s_2\,s_4\,(1+s_3+s_2\,s_3\,s_4)}
{(1+s_2\,s_4)\,(1+s_2\,s_3\,s_4)}\,,\quad
N=\frac1{s_2^2\,s_4}\,,\quad
\rho=2^{-{1/2}}\,s_2^{1/2}\,(1+s_2\,s_4),\quad
\si=-\,\frac{1+s_3-s_2^2\,s_4}{s_2\,(1+s_2\,s_4)}\,,\\
&\rho_0=\frac{(1+s_3+s_2\,s_3\,s_4)^{1/2}}{1+s_2\,s_3\,s_4}\,,\quad
\si_0=-\,\frac{s_2^{1/2}\,s_3\,s_4\,(s_2+s_3)}
{2^{1/2}\,(1+s_3+s_2\,s_3\,s_4)^{1/2}}\,.
\end{align*}
{\bf Charlier}:
\begin{align*}
&p_n(x;0,s_2,s_3,0)
=\rho_0^n\,p_n\left(\rho_0^{-1}x;0,\frac{s_2}{1+s_3},0,0\right)=
\rho^n\,c_n(\rho^{-1}x-\si;a),\sLP
&a=\frac{1+s_3}{s_2}\,,\quad
\rho=2^{-{1/2}}\,s_2^{1/2}\,,\quad
\si=-\,\frac{1+s_3}{s_2}\,,\quad
\rho_0=(1+s_3)^{1/2}.
\end{align*}
{\bf Hermite}:
\begin{align*}
&p_n(x;s_1,0,s_3,s_4)=p_n(x;s_1,0,s_3,0)=p_n(x;0,0,s_3,0)=
\rho_0^n\,p_n(\rho_0^{-1}\,x;0,0,s_3,s_4)\\
&=\rho_0^n\,p_n(\rho_0^{-1}\,x;0,0,s_3,0)
=\rho^n\,p_n(\rho^{-1}\,x;s_1,0,0,s_4)=\rho^n\,p_n(\rho^{-1}\,x;0,0,0,s_4)=
\rho^n\,h_n(\rho^{-1}\,x),\\
&\rho=\frac{(1+s_3)^{1/2}\,(1+s_1+s_1\,s_3)^{1/2}}{(1+s_1)^{1/2}}\,,\quad
\rho_0=\frac{(1+s_1+s_1\,s_3)^{1/2}}{(1+s_1)^{1/2}}\,.
\end{align*}

\begin{figure}[t]
\hskip1cm
\setlength{\unitlength}{3mm}
\begin{picture}(0,37.5)
\put(22,34.5) {\framebox(5,3) {Racah}}
\put(24.5,34.5) {\color{red}\vector(-2,-1){10}}
\put(16.2,32){\color{red}{$s_1\downarrow0$}}
\put(24.5,34.5) {\color{blue}\vector(-1,-2){2.5}}
\put(24.5,34.5) {\color{blue}\vector(-1,-2){2.2}}
\put(21.2,32){\color{blue}{$s_2\downarrow0$}}
\put(24.5,34.5) {\color{cyan}\vector(1,-1){5}}
\put(24.5,34.5) {\color{cyan}\vector(1,-1){4.5}}
\put(24.5,34.5) {\color{cyan}\vector(1,-1){4}}
\put(25.7,32){\color{cyan}{$s_3\downarrow0$}}
\put(24.5,34.5) {\color{black}\vector(3,-1){15}}
\put(24.5,34.5) {\color{black}\vector(3,-1){14}}
\put(24.5,34.5) {\color{black}\vector(3,-1){13}}
\put(24.5,34.5) {\color{black}\vector(3,-1){12}}
\put(29.6,32){\color{black}{$s_4\downarrow0$}}
\put(10,26.5) {\framebox(7.5,3) {Krawtchouk}}
\put(20,26.5) {\framebox(5.5,3) {Hermite}}
\put(27,26.5) {\framebox(7,3) {dual Hahn}}
\put(36,26.5) {\framebox(5,3) {Meixner}}
\put(12.5,26.5) {\color{cyan}\vector(-3,-2){7.5}}
\put(12.5,26.5) {\color{cyan}\vector(-3,-2){6.9}}
\put(12.5,26.5) {\color{cyan}\vector(-3,-2){6.3}}
\put(13.5,26.5) {\color{blue}\vector(0,-1){5}}
\put(13.5,26.5) {\color{blue}\vector(0,-1){4.5}}
\put(14.5,26.5) {\color{black}\vector(14,-5){14}}
\put(14.5,26.5) {\color{black}\vector(14,-5){13}}
\put(14.5,26.5) {\color{black}\vector(14,-5){12}}
\put(14.5,26.5) {\color{black}\vector(14,-5){11}}
\put(23.5,26.5) {\color{red}\vector(-9,-5){9}}
\put(23.5,26.5) {\color{cyan}\vector(-1,-2){2.5}}
\put(23.5,26.5) {\color{cyan}\vector(-1,-2){2.2}}
\put(23.5,26.5) {\color{cyan}\vector(-1,-2){1.9}}
\put(23.5,26.5) {\color{black}\vector(3,-1){15}}
\put(23.5,26.5) {\color{black}\vector(3,-1){14}}
\put(23.5,26.5) {\color{black}\vector(3,-1){13}}
\put(23.5,26.5) {\color{black}\vector(3,-1){12}}
\put(28,26.5) {\color{red}\vector(-4,-1){20}}
\put(28,26.5) {\color{blue}\vector(-1,-1){5}}
\put(28,26.5) {\color{blue}\vector(-1,-1){4.5}}
\put(29,26.5) {\color{black}\vector(3,-1){15}}
\put(29,26.5) {\color{black}\vector(3,-1){14}}
\put(29,26.5) {\color{black}\vector(3,-1){13}}
\put(29,26.5) {\color{black}\vector(3,-1){12}}
\put(39,26.5) {\color{red}\vector(-2,-1){10}}
\put(39,26.5) {\color{blue}\vector(0,-1){5}}
\put(39,26.5) {\color{blue}\vector(0,-1){4.5}}
\put(39,26.5) {\color{cyan}\vector(1,-1){5}}
\put(39,26.5) {\color{cyan}\vector(1,-1){4.5}}
\put(39,26.5) {\color{cyan}\vector(1,-1){4}}
\put(2,18.5) {\framebox(7.5,3) {Krawtchouk}}
\put(11,18.5) {\framebox(5.5,3) {Hermite}}
\put(18,18.5) {\framebox(5.5,3) {Hermite}}
\put(26,18.5) {\framebox(5.5,3) {Charlier}}
\put(34,18.5) {\framebox(5.5,3) {Hermite}}
\put(41.5,18.5) {\framebox(5,3) {Meixner}}
\put(3.7,18.5) {\color{blue}{\vector(2,-1){10}}}
\put(3.7,18.5) {\color{blue}{\vector(2,-1){9.3}}}
\put(5,18.5) {\color{black}{\vector(3,-1){15}}}
\put(5,18.5) {\color{black}{\vector(3,-1){14}}}
\put(5,18.5) {\color{black}{\vector(3,-1){13}}}
\put(5,18.5) {\color{black}{\vector(3,-1){12}}}
\put(14,18.5) {\color{cyan}\vector(0,-1){5}}
\put(14,18.5) {\color{cyan}\vector(0,-1){4.5}}
\put(14,18.5) {\color{cyan}\vector(0,-1){4}}
\put(14,18.5) {\color{black}\vector(3,-1){15}}
\put(14,18.5) {\color{black}\vector(3,-1){14}}
\put(14,18.5) {\color{black}\vector(3,-1){13}}
\put(14,18.5) {\color{black}\vector(3,-1){12}}
\put(19.5,18.5) {\color{red}\vector(-1,-1){5}}
\put(21.6,18.5) {\color{black}\vector(3,-1){15}}
\put(21.6,18.5) {\color{black}\vector(3,-1){14}}
\put(21.6,18.5) {\color{black}\vector(3,-1){13}}
\put(21.6,18.5) {\color{black}\vector(3,-1){12}}
\put(28,18.5) {\color{cyan}\vector(-1,-1){5}}
\put(28,18.5) {\color{cyan}\vector(-1,-1){4.5}}
\put(28,18.5) {\color{cyan}\vector(-1,-1){4}}
\put(29.5,18.5) {\color{blue}\vector(0,-1){5}}
\put(29.5,18.5) {\color{blue}\vector(0,-1){4.5}}
\put(35.4,18.5) {\color{red}\vector(-1,-1){5}}
\put(37.8,18.5) {\color{cyan}\vector(0,-1){5}}
\put(37.8,18.5) {\color{cyan}\vector(0,-1){4.5}}
\put(37.8,18.5) {\color{cyan}\vector(0,-1){4}}
\put(44,18.5) {\color{red}\vector(-4,-1){20}}
\put(44,18.5) {\color{blue}\vector(-1,-1){5}}
\put(44,18.5) {\color{blue}\vector(-1,-1){4.5}}
\put(11.5,10.5) {\framebox(5.5,3) {Hermite}}
\put(19,10.5) {\framebox(5.5,3) {Charlier}}
\put(27,10.5) {\framebox(5.5,3) {Hermite}}
\put(35,10.5) {\framebox(5.5,3) {Hermite}}
\put(15,10.5) {\color{black}\vector(2,-1){10}}
\put(15,10.5) {\color{black}\vector(2,-1){9.3}}
\put(15,10.5) {\color{black}\vector(2,-1){8.6}}
\put(15,10.5) {\color{black}\vector(2,-1){7.9}}
\put(23,10.5) {\color{blue}\vector(1,-2){2.5}}
\put(23,10.5) {\color{blue}\vector(1,-2){2.2}}
\put(28.5,10.5) {\color{cyan}\vector(-1,-2){2.5}}
\put(28.5,10.5) {\color{cyan}\vector(-1,-2){2.2}}
\put(28.5,10.5) {\color{cyan}\vector(-1,-2){1.9}}
\put(37,10.5) {\color{red}\vector(-2,-1){10}}
\put(23,2.5) {\framebox(5.5,3) {Hermite}}
\end{picture}
\vskip-1cm
\caption{The various parameter restrictions in the second chart for
the Racah manifold}
\label{fig:5}
\end{figure}
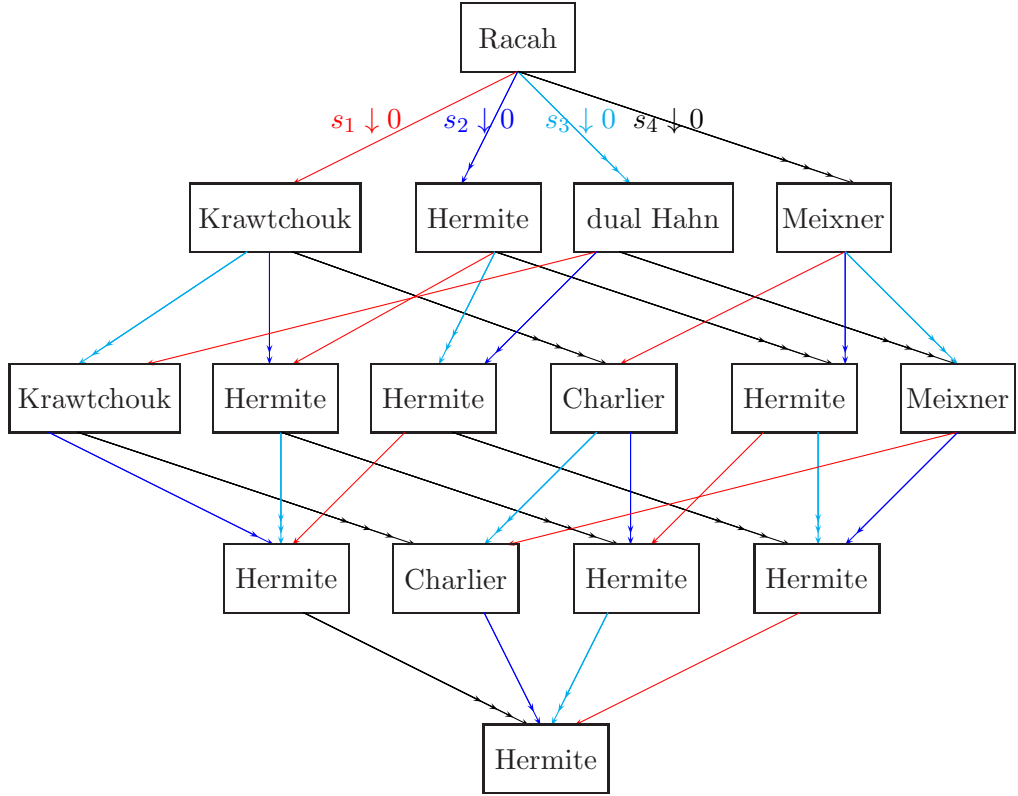

The various parameter restrictions of the polynomial \eqref{37}
are summarized in Figure \ref{fig:5}. Even more than in
Figure \ref{fig:4}, several boxes essentially coincide.

\subsection{From Racah to Hermite along dual Hahn and Laguerre}
Here we will see the chart for the Racah manifold corresponding to the
third graph in Figure \ref{fig:2}.
Let $r_n$ again be the monic Racah polynomial given by \eqref{27}.
For $u_1,u_2,u_3,u_4>0$, $u_2^2\,u_3\,u_4<1$ and
$u_2\,u_3\,(u_1-u_2)<1$ put
\begin{align}
&p_n(x)=p_n(x;u_1,u_2,u_3,u_4):=
\rho^n\,r_n(\rho^{-1}x-\si;\al,\be,-N-1,\de),
\label{38}\sLP
{\rm where}\quad
&\al=\frac{1+u_1}{u_2}\,,\quad
\be=\frac1{u_1\,u_2^2\,u_3^2\,u_4}\,,\quad
N=\frac1{u_2^2\,u_3\,u_4}\,,\quad
\de=\frac{1+u_4+u_2\,u_3\,u_4+u_2^2\,u_3\,u_4}{u_2^2\,u_3\,u_4}\,,
\label{41}\sLP
&\rho=2^{-{1/2}}\,u_2^{5/2}\,u_3\,u_4\,,\quad
\si=-\,\frac{(1+u_1)\,(1+u_1\,u_3+u_1\,u_3\,u_4)}{u_2^3\,u_3\,u_4}\,.
\nonu
\end{align}
Then the inequalities for $t_1,t_2,t_3,t_4$ are equivalent with \eqref{18}.
Furthermore \eqref{32} holds with
\begin{align*}
&B_n=-2^{-{1/2}}\,u_2^{1/2}\,
\Big(2\,n^2\,(n+1)^2\,u_1^2\,u_2^6\,u_3^5\,u_4^3
+4\,n^2\,(n+1)\,u_1\,u_2^4\,u_3^3\,u_4^2\,(1+u_1\,u_2\,u_3^2\,u_4
+u_1^2\,u_2\,u_3^2\,u_4)\\
&\qquad\qquad+n^2\,u_2^2\,u_3\,u_4\,
\big(2-2\,u_1\,u_3-2\,u_1^2\,u_3^2-u_1\,u_3\,u_4-2\,u_1^2\,u_3^2\,u_4
+5\,u_1\,u_2\,u_3^2\,u_4\\
&\qquad\qquad\qquad+6\,u_1^2\,u_2\,u_3^3\,u_4
+2\,u_1^2\,u_2\,u_3^3\,u_4+4u_1^3\,u_2\,u_3^3\,u_4
-4\,u_1^2\,u_2^2\,u_3^3\,u_4+4\,u_1^3\,u_2\,u_3^4\,u_4\\
&\qquad\qquad\qquad+4\,u_1^4\,u_2\,u_3^4\,u_4+u_1^2\,u_2\,u_3^3\,u_4^2
+u_1^3\,u_2\,u_3^3\,u_4^2+4\,u_1^3\,u_2\,u_3^4\,u_4^2
+4\,u_1^4\,u_2\,u_3^4\,u_4^2\\
&\qquad\qquad\qquad+3\,u_1^2\,u_2^2\,u_3^4\,u_4^2
+5\,u_1^3\,u_2^2\,u_3^4\,u_4^2+2\,u_1^4\,u_2^2\,u_3^4\,u_4^2
+2\,u_1^2\,u_2^3\,u_3^4\,u_4^2+2\,u_1^3\,u_2^3\,u_3^4\,u_4^2
\big)\displaybreak[0]\\
&\qquad\qquad+n\,(1+u_1\,u_2\,u_3^2\,u_4+u_1^2\,u_2\,u_3^2\,u_4
+u_1\,u_2^2\,u_3^2\,u_4)\,\big(
-2 - 2\,u_1\,u_3 - u_4 - 2\,u_1\,u_3\,u_4\\
&\qquad\qquad\qquad + u_2\,u_3\,u_4 +
2\,u_1\,u_2\,u_3\,u_4 +2\,u_1\,u_2\,u_3^2\,u_4+4\,u_1^2\,u_2\,u_3^2\,u_4
-4\,u_1\,u_2^2\,u_3^2\,u_4\\
&\qquad\qquad\qquad
+4\,u_1^2\,u_2\,u_3^3\,u_4+4\,u_1^3\,u_2\,u_3^3\,u_4
+u_1\,u_2\,u_3^2\,u_4^2 + 
 u_1^2\,u_2\,u_3^2\,u_4^2 + 4\,u_1^2\,u_2\,u_3^3\,u_4^2\\
&\qquad\qquad\qquad+ 4\,u_1^3\,u_2\,u_3^3\,u_4^2 + 
u_1\,u_2^2\,u_3^3\,u_4^2 + u_1^2\,u_2^2\,u_3^3\,u_4^2 +
2\,u_1\,u_2^3\,u_3^3\,u_4^2 + 
2\,u_1^2\,u_2^3\,u_3^3\,u_4^2\big)
\displaybreak[0]\\
&\qquad\qquad
+(1 + u_1\,u_2\,u_3^2\,u_4 + u_1^2\,u_2\,u_3^2\,u_4)\,
\big(-1 - u_1\,u_3 - u_1\,u_3\,u_4 + 
u_1^2\,u_3^2\,u_4 + u_1^3\,u_3^2\,u_4\\
&\qquad\qquad\qquad 
-u_1\,u_2\,u_3^2\,u_4 - 2\,u_1\,u_2^2\,u_3^2\,u_4
+u_1^2\,u_3^3\,u_4 + 2\,u_1^3\,u_3^3\,u_4 + u_1^4\,u_3^3\,u_4
+ 2\,u_1^2\,u_2\,u_3^3\,u_4\\
&\qquad\qquad\qquad 
+2\,u_1^3\,u_2\,u_3^3\,u_4 + u_1^2\,u_3^3\,u_4^2
+ 2\,u_1^3\,u_3^3\,u_4^2
+u_1^4\,u_3^3\,u_4^2 + 2\,u_1^2\,u_2\,u_3^3\,u_4^2 +
2\,u_1^3\,u_2\,u_3^3\,u_4^2\big)
\Big)
\displaybreak[0]\\
&\qquad\quad\times\big(1+u_1\,u_2\,u_3^2\,u_4\,(1+u_1)
+2\,n\,u_1\,u_2^2\,u_3^2\,u_4\big)^{-1}\\
&\qquad\quad\times\big(1+u_1\,u_2\,u_3^2\,u_4\,(1+u_1)
+2\,(n+1)\,u_1\,u_2^2\,u_3^2\,u_4\big)^{-1},
\displaybreak[0]\mLP
&C_n=\thalf\,n\,(1+u_1+n\,u_2)\,(1+(1-n)\,u_2^2\,u_3\,u_4)\,
(1+u_4-u_1\,u_2\,u_3\,u_4+(1-n)\,u_2^2\,u_3\,u_4)\sLP
&\quad\times\frac{
\big(1+u_1\,u_2\,u_3^2\,u_4\,(1+u_1)+n\,u_1\,u_2^2\,u_3^2\,u_4\big)\,
\big(1+u_1\,u_3\,(1+u_4+u_2\,u_3\,u_4)+(n+1)\,u_1\,u_2^2\,u_3^2\,u_4\big)}
{\big(1+u_1\,u_2\,u_3^2\,u_4\,(1+u_1)
+(2\,n-1)\,u_1\,u_2^2\,u_3^2\,u_4\big)\,
\big(1+u_1\,u_2\,u_3^2\,u_4\,(1+u_1)
+2\,n\,u_1\,u_2^2\,u_3^2\,u_4\big)^2}
\displaybreak[0]\sLP
&\quad\times\frac{(1+n\,u_1\,u_2^2\,u_3^2\,u_4)\,
\big(1+u_1\,u_3\,(1+u_2\,u_3\,u_4+u_1\,u_2\,u_3\,u_4)
+(n+1)\,u_1\,u_2^2\,u_3^2\,u_4\big)}
{1+u_1\,u_2\,u_3^2\,u_4\,(1+u_1)
+(2\,n+1)\,u_1\,u_2^2\,u_3^2\,u_4}\,.
\end{align*}
Note that $B_n$ and $C_n$, as functions of $u_1,u_2,u_3,u_4>0$,
can be uniquely extended to continuous functions of
$u_1,u_2,u_3,u_4\ge0$.

\begin{figure}[t]
\setlength{\unitlength}{3mm}
\begin{picture}(0,37.5)
\put(22,34.5) {\framebox(5,3) {Racah}}
\put(24.5,34.5) {\color{red}\vector(-2,-1){10}}
\put(16.2,32){\color{red}{$u_1\downarrow0$}}
\put(24.5,34.5) {\color{blue}\vector(-1,-2){2.5}}
\put(24.5,34.5) {\color{blue}\vector(-1,-2){2.2}}
\put(21,32){\color{blue}{$u_2\downarrow0$}}
\put(24.5,34.5) {\color{cyan}\vector(1,-1){5}}
\put(24.5,34.5) {\color{cyan}\vector(1,-1){4.5}}
\put(24.5,34.5) {\color{cyan}\vector(1,-1){4}}
\put(25.5,32){\color{cyan}{$u_3\downarrow0$}}
\put(24.5,34.5) {\color{black}\vector(3,-1){15}}
\put(24.5,34.5) {\color{black}\vector(3,-1){14}}
\put(24.5,34.5) {\color{black}\vector(3,-1){13}}
\put(24.5,34.5) {\color{black}\vector(3,-1){12}}
\put(29.4,32){\color{black}{$u_4\downarrow0$}}
\put(10,26.5) {\framebox(7,3) {dual Hahn}}
\put(20,26.5) {\framebox(5.5,3) {Hermite}}
\put(27,26.5) {\framebox(6,3) {Meixner}}
\put(35,26.5) {\framebox(6.5,3) {Laguerre}}
\put(13.5,26.5) {\color{cyan}\vector(-2,-1){10}}
\put(13.5,26.5) {\color{cyan}\vector(-2,-1){9.3}}
\put(13.5,26.5) {\color{cyan}\vector(-2,-1){8.6}}
\put(13.5,26.5) {\color{blue}\vector(0,-1){5}}
\put(13.5,26.5) {\color{blue}\vector(0,-1){4.5}}
\put(13.5,26.5) {\color{black}\vector(3,-1){15}}
\put(13.5,26.5) {\color{black}\vector(3,-1){14}}
\put(13.5,26.5) {\color{black}\vector(3,-1){13}}
\put(13.5,26.5) {\color{black}\vector(3,-1){12}}
\put(23.5,26.5) {\color{red}\vector(-2,-1){10}}
\put(23.5,26.5) {\color{cyan}\vector(-1,-2){2.5}}
\put(23.5,26.5) {\color{cyan}\vector(-1,-2){2.2}}
\put(23.5,26.5) {\color{cyan}\vector(-1,-2){1.9}}
\put(23.5,26.5) {\color{black}\vector(3,-1){15}}
\put(23.5,26.5) {\color{black}\vector(3,-1){14}}
\put(23.5,26.5) {\color{black}\vector(3,-1){13}}
\put(23.5,26.5) {\color{black}\vector(3,-1){12}}
\put(28,26.5) {\color{red}\vector(-4,-1){20}}
\put(28,26.5) {\color{blue}\vector(-1,-1){5}}
\put(28,26.5) {\color{blue}\vector(-1,-1){4.5}}
\put(29,26.5) {\color{black}\vector(3,-1){15}}
\put(29,26.5) {\color{black}\vector(3,-1){14}}
\put(29,26.5) {\color{black}\vector(3,-1){13}}
\put(29,26.5) {\color{black}\vector(3,-1){12}}
\put(39,26.5) {\color{red}\vector(-2,-1){10}}
\put(39,26.5) {\color{blue}\vector(0,-1){5}}
\put(39,26.5) {\color{blue}\vector(0,-1){4.5}}
\put(39.5,26.5) {\color{cyan}\vector(1,-1){5}}
\put(39.5,26.5) {\color{cyan}\vector(1,-1){4.5}}
\put(39.5,26.5) {\color{cyan}\vector(1,-1){4}}
\put(2.5,18.5) {\framebox(6,3) {Meixner}}
\put(10,18.5) {\framebox(5.5,3) {Hermite}}
\put(18,18.5) {\framebox(5.5,3) {Hermite}}
\put(26,18.5) {\framebox(6.5,3) {Laguerre}}
\put(34,18.5) {\framebox(6,3) {Hermite}}
\put(41,18.5) {\framebox(6.5,3) {Laguerre}}
\put(3.5,18.5) {\color{blue}{\vector(2,-1){10}}}
\put(3.5,18.5) {\color{blue}{\vector(2,-1){9.3}}}
\put(6,18.5) {\color{black}{\vector(3,-1){15}}}
\put(6,18.5) {\color{black}{\vector(3,-1){14}}}
\put(6,18.5) {\color{black}{\vector(3,-1){13}}}
\put(6,18.5) {\color{black}{\vector(3,-1){12}}}
\put(14,18.5) {\color{cyan}\vector(0,-1){5}}
\put(14,18.5) {\color{cyan}\vector(0,-1){4.5}}
\put(14,18.5) {\color{cyan}\vector(0,-1){4}}
\put(14,18.5) {\color{black}\vector(3,-1){15}}
\put(14,18.5) {\color{black}\vector(3,-1){14}}
\put(14,18.5) {\color{black}\vector(3,-1){13}}
\put(14,18.5) {\color{black}\vector(3,-1){12}}
\put(19.5,18.5) {\color{red}\vector(-1,-1){5}}
\put(21.6,18.5) {\color{black}\vector(3,-1){15}}
\put(21.6,18.5) {\color{black}\vector(3,-1){14}}
\put(21.6,18.5) {\color{black}\vector(3,-1){13}}
\put(21.6,18.5) {\color{black}\vector(3,-1){12}}
\put(28,18.5) {\color{cyan}\vector(-1,-1){5}}
\put(28,18.5) {\color{cyan}\vector(-1,-1){4.5}}
\put(28,18.5) {\color{cyan}\vector(-1,-1){4}}
\put(29.2,18.5) {\color{blue}\vector(0,-1){5}}
\put(29.2,18.5) {\color{blue}\vector(0,-1){4.5}}
\put(34.4,18.5) {\color{red}\vector(-1,-1){5}}
\put(37,18.5) {\color{cyan}\vector(0,-1){5}}
\put(37,18.5) {\color{cyan}\vector(0,-1){4.5}}
\put(37,18.5) {\color{cyan}\vector(0,-1){4}}
\put(44,18.5) {\color{red}\vector(-4,-1){20}}
\put(45,18.5) {\color{blue}\vector(-1,-1){5}}
\put(45,18.5) {\color{blue}\vector(-1,-1){4.5}}
\put(12.5,10.5) {\framebox(5.5,3) {Hermite}}
\put(19,10.5) {\framebox(6.5,3) {Laguerre}}
\put(27,10.5) {\framebox(5.5,3) {Hermite}}
\put(35,10.5) {\framebox(6,3) {Hermite}}
\put(15,10.5) {\color{black}\vector(2,-1){10}}
\put(15,10.5) {\color{black}\vector(2,-1){9.3}}
\put(15,10.5) {\color{black}\vector(2,-1){8.6}}
\put(15,10.5) {\color{black}\vector(2,-1){7.9}}
\put(23,10.5) {\color{blue}\vector(1,-2){2.5}}
\put(23,10.5) {\color{blue}\vector(1,-2){2.2}}
\put(28.5,10.5) {\color{cyan}\vector(-1,-2){2.5}}
\put(28.5,10.5) {\color{cyan}\vector(-1,-2){2.2}}
\put(28.5,10.5) {\color{cyan}\vector(-1,-2){1.9}}
\put(37,10.5) {\color{red}\vector(-2,-1){10}}
\put(23.0,2.5) {\framebox(5.5,3) {Hermite}}
\end{picture}
\vskip-1cm
\caption{The various parameter restrictions in the third chart for the
Racah manifold}
\label{fig:6}
\end{figure}
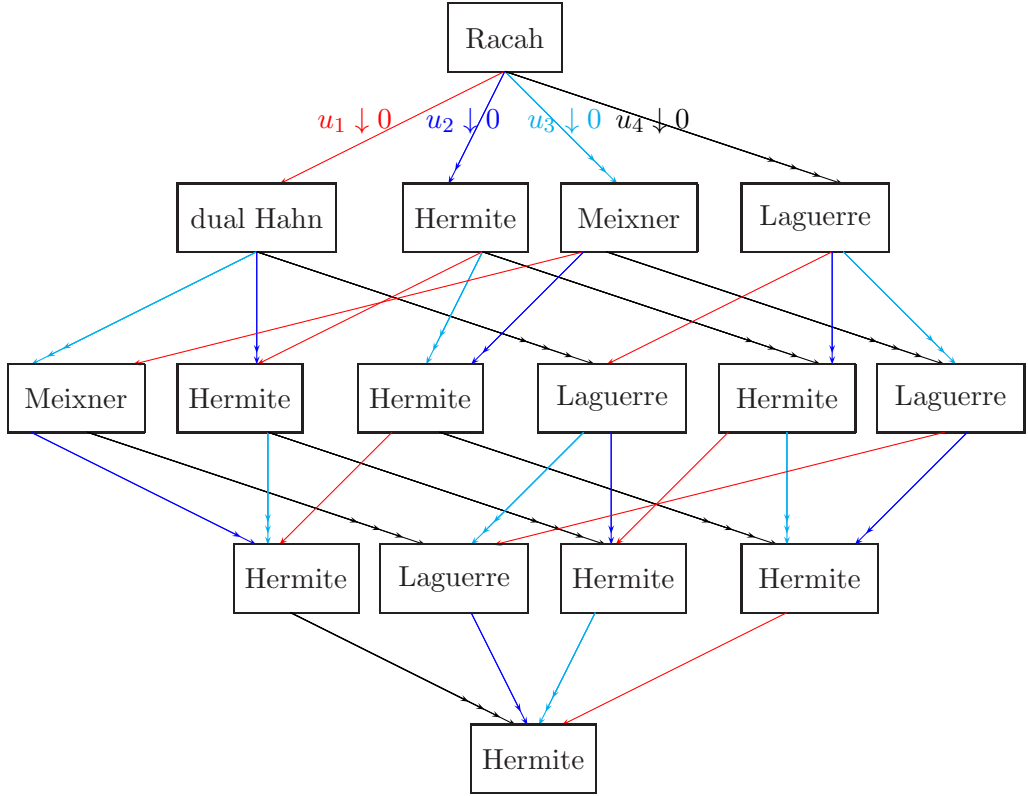

We now put one or more of the $u_1,u_2,u_3,u_4$ equal to zero and
we proceed as in \S\ref{35}. We obtain:
\mLP
{\bf dual Hahn}:
\begin{align*}
&p_n(x;0,u_2,u_3,u_4)=
\rho^n\,r_n^{\rm DH}(\rho^{-1}x-\si;\ga,\de,N),\sLP
&\ga=\frac1{u_2}\,,\quad
\de=\frac1{u_2^2\,u_3}\,,\quad
N=\frac1{u_2^2\,u_3\,u_4}\,,\quad
\rho=2^{-{1/2}}\,u_2^{5/2}\,u_3\,u_4,\quad
\si=-\,\frac1{u_2^2\,u_3\,u_4}\,.
\end{align*}
{\bf Meixner}:
\begin{align*}
&p_n(x;u_1,u_2,0,u_4)=\rho_0^n\,
p_n\Big(\rho_0^{-1}x;0,\frac{u_2}{1+u_1}\,,0,u_4\Big)
=\rho^n\,m_n(\rho^{-1}x-\si;\be,c),\sLP
&\be=\frac{1+u_1+u_2}{u_2}\,,\quad
c=\frac1{1+u_4},\quad
\rho=2^{-{1/2}}\,u_2^{1/2}\,u_4,\quad
\rho_0=(1+u_1)^{1/2},\quad
\si=-\,\frac{1+u_1}{u_2\,u_4}\,.
\end{align*}
{\bf Laguerre}:
\begin{align*}
&p_n(x;u_1,u_2,u_3,0)
=\rho_1^n\,p_n(\rho_1^{-1}x;u_1,u_2,0,0)
=\rho_0^n\,p_n(\rho_0^{-1}x;0,\frac{u_2}{1+u_1}\,,u_3,0)\\
&\qquad\qquad=\rho_0^n\,p_n(\rho_0^{-1}x;0,\frac{u_2}{1+u_1}\,,0,0)
=\rho^n\,\ell_n^{(\al)}(\rho^{-1}x-\si),\sLP
&\al=\frac{1+u_1}{u_2}\,,\quad
\rho=2^{-{1/2}}\,u_2^{1/2}\,(1+u_1\,u_3),\quad
\rho_0=(1+u_1)^{1/2}\,(1+u_1\,u_3),\\
&\rho_1=1+u_1\,u_3,\quad
\si=-\,\frac{1+u_1}{u_2}\,.
\end{align*}
{\bf Hermite}:
\begin{align*}
&p_n(x;u_1,0,u_3,u_4)
=\rho_0^n\,p_n(\rho_0^{-1}\,x;0,0,u_3,u_4)
=\rho_0^n\,p_n(\rho_0^{-1}\,x;0,0,0,u_4)\\
&\qquad\qquad=\rho_1^n\,p_n(\rho_1^{-1}\,x;u_1,0,0,u_4)
=\rho_2^n\,p_n(\rho_2^{-1}\,x;u_1,0,u_3,0)
=\rho_3^n\,p_n(\rho_3^{-1}\,x;u_1,0,0,0)\\
&\qquad\qquad=\rho^n\,p_n(\rho^{-1}x;0,0,u_3,0)
=\rho^n\,p_n(\rho^{-1}x;0,0,0,0)
=\rho^n\,h_n(\rho^{-1}x),
\displaybreak[0]\sLP
&\rho=(1+u_1)^{1/2}\,(1+u_1\,u_3)^{1/2}\,(1+u_4)^{1/2}\,
\big(1+u_1\,u_3\,(1+u_4)\big)^{1/2},\\
&\rho_0=
(1+u_1)^{1/2}\,(1+u_1\,u_3)^{1/2}\,\big(1+u_1\,u_3\,(1+u_4)\big)^{1/2},\;\;
\rho_1=(1+u_1\,u_3)^{1/2}\,\big(1+u_1\,u_3\,(1+u_4)\big)^{1/2},\\
&\rho_2=\frac{(1+u_4)^{1/2}\,
\big(1+u_1\,u_3\,(1+u_4)\big)^{1/2}}{(1+u_1\,u_3)^{1/2}}\,,\;\;
\rho_3=
(1+u_4)^{1/2}\,(1+u_1\,u_3)^{1/2}\,\big(1+u_1\,u_3\,(1+u_4)\big)^{1/2}.
\end{align*}
The various parameter restrictions of the polynomial \eqref{38}
are summarized in Figure \ref{fig:6}. Many boxes essentially coincide.

\subsection{Transformations between the local charts}
By \eqref{39}, \eqref{40} and \eqref{41} we can compute the
transformations between the local coordinates $(t_1,t_2,t_3,t_4)$,
$(s_1,s_2,s_3,s_4)$ and $(u_1,u_2,u_3,u_4)$.\sLP
{\bf first chart $\leftrightarrow$ second chart}:
\begin{align*}
&t_1=\frac{s_1\,s_2}{1+s_1}\,,\quad
t_2=s_2\,s_3\,s_4,\quad
t_3=\frac{1+s_1}{s_2(1+s_1+s_1\,s_2^2\,s_4)}\,,\quad
t_4=\frac{s_2}{s_3}\,;\\
&s_1=\frac{t_1\,t_3}{1-t_1\,t_3-t_1\,t_2\,t_3\,t_4}\,,\;\;\;
s_2=\frac{1-t_1\,t_2\,t_3\,t_4}{t_3}\,,\;\;\;
s_3=\frac{1-t_1\,t_2\,t_3\,t_4}{t_3\,t_4}\,,\;\;\;
s_4=\frac{t_2\,t_3^2\,t_4}{(1-t_1\,t_2\,t_3\,t_4)^2}\,.
\end{align*}
This is a homeomorphism between
\[
\{(t_1,t_2,t_3,t_4)\mid t_1\ge0,\;t_2\ge0,\;t_3>0,\;t_4>0,\;
t_1t_3(1+t_2t_4)<1,\;t_2\,t_4<1\}
\]
and
\[
\{(s_1,s_2,s_3,s_4)\mid s_1\ge0,\;s_2>0,s_3>0,\;s_4\ge0,\;s_2^2\,s_4<1\}.
\]
It also identifies the Krawtchouk box $\{t_1=0\}$, Meixner box
$\{t_2=0\}$ and Charlier box $\{t_1=t_2=0\}$ in
Figure \ref{fig:4} with the corresponding boxes
$\{s_1=0\}$, $\{s_4=0\}$ and $\{s_1=s_4=0\}$ in Figure \ref{fig:5}.
\mLP
{\bf second chart $\leftrightarrow$ third chart}:
\begin{align*}
&s_1=\frac1{u_4\,(1-u_1\,u_2\,u_3)}\,,\quad
s_2=\frac{u_2\,\big(1+u_4\,(1-u_1\,u_2\,u_3)\big)}{1+u_1}\,,\\
&s_3=u_1\,u_3\,\big(1+u_4\,(1-u_1\,u_2\,u_3)\big),\quad
s_4=\frac{(1+u_1)^2\,u_3\,u_4}{\big(1+u_4\,(1-u_1\,u_2\,u_3)\big)^2}\,;\\
&u_1=\frac{-2s_1^2\,s_3+s_4\,(1+s_1)\,(1 + 2\,s_1 + s_1^2 - s_1^2\,s_2\,s_3)
-(1+s_1)\,S^{1/2}}
{2\,s_1^2\,s_3\,(1 + s_2\,s_4 + s_1\,s_2\,s_4)}\,,\\
&u_2=\frac{s_1\,s_2^2\,s_4}{1+s_2\,s_4+s_1\,s_2\,s_4}+
s_2\,\frac{s_4\,(1 + 2\,s_1 + s_1^2 - s_1^2\,s_2\,s_3)-S^{1/2}}
{2\,s_1\,s_3\,(1 + s_2\,s_4 + s_1\,s_2\,s_4)}\,,\\
&u_3=\frac{-2s_1^2\,s_3+s_4\,(1+s_1)\,(1 + 2\,s_1 + s_1^2 - s_1^2\,s_2\,s_3)
+(1+s_1)\,S^{1/2}}
{2\,s_1\,(1+s_1)}\,,\\
&u_4=\frac1{s_1}+s_2\,\frac{s_4\,
(1 + 2\,s_1 + s_1^2 - s_1^2\,s_2\,s_3)-S^{1/2}}
{2\,s_1\,(1+s_1)}\,,\\
\intertext{where}
&S:=s_4^2\,(1 + 2\,s_1 + s_1^2 - s_1^2\,s_2\,s_3)^2 -
4\,s_1^2\,s_3\,s_4\,(1 + s_1).
\end{align*}
This is a homeomorphism between
\begin{multline*}
\{(s_1,s_2,s_3,s_4)\mid s_1>0,\;s_2\ge0,s_3>0,\;s_4>0,\;s_2^2\,s_4<1,\;
S\ge0,\\
-2\,s_1^2\,s_3+s_4\,(1+s_1)\,(1+2\,s_1+s_1^2-s_1^2\,s_2\,s_3)\ge0\,\}
\end{multline*}
and
\begin{multline*}
\{u_1,u_2,u_3,u_4>0\mid
u_1>0,\;u_2\ge0,\;u_3>0,\;u_4>0,\;u_2^2\,u_3\,u_4<1,\\
u_2\,u_3\,(u_1-u_2)<1,\;u_1\,(1+u_2\,u_3)\le1\}.
\end{multline*}
It also identifies the Hermite box $\{s_2=0\}$ in Figure \ref{fig:5}
with the Hermite box $\{u_2=0\}$ in Figure~\ref{fig:6}.
\mLP
{\bf first chart $\leftrightarrow$ third chart}:
\begin{align*}
&t_1=\frac{u_2}{1+u_1}\,,\quad
t_2=u_1\,u_2\,u_3^2\,u_4\,(1 + u_1),\\
&t_3=\frac{1+u_1}{u_2\,\big(1 + u_4\,(1-u_2\,u_3\,(u_1-u_2))\big)}\,,
\quad
t_4=\frac{u_2}{u_1\,u_3\,(1+u_1)}\,;\displaybreak[0]\\
&u_1=\frac{t_2\,t_3\,t_4\,(t_4-t_1^2)
-2\,t_1^2\,\big(1-t_1\,t_3\,(1+t_2\,t_4)\big)- t_4\,T^{1/2}}
{2\,t_1^2\,\big(1-t_1\,t_3\,(1+t_2\,t_4)+t_2\,t_3\,t_4\big)}\,,
\displaybreak[0]\\
&u_2=\frac{t_2\,t_3\,t_4\,(t_1^2 + t_4) - t_4\,T^{1/2}}
{2\,t_1\,\big(1-t_1\,t_3\,(1+t_2\,t_4)+t_2\,t_3\,t_4\big)}\,,
\displaybreak[0]\\
&u_3=\frac{t_2\,t_3\,t_4\,(t_4-t_1^2)
-2\,t_1^2\,\big(1-t_1\,t_3\,(1+t_2\,t_4)\big)+t_4\,T^{1/2}}
{2\,t_1\,\big(1 - t_1\,t_3\,(1+t_2\,t_4)\big)}\,,\\
&u_4=\frac{2\big(1-t_1\,t_3\,(1+t_2\,t_4)\big)+t_2\,t_3\,(t_4-t_1^2)-T^{1/2}}
{2\,t_1\,t_3}\,,\\
\intertext{where}
&T:=t_2^2\,t_3^2\,(t_4-t_1^2)^2-4\,t_1^2\,t_2\,t_3
\big(1-t_1\,t_3\,(1+t_2\,t_4)\big).
\end{align*}
This is a homeomorphism between
\begin{multline*}
\{(t_1,t_2,t_3,t_4)\mid t_1>0,\;t_2>0,\;t_3>0,\;t_4>0,\;t_2\,t_4<1,\;T\ge0,\\
t_2\,t_3\,t_4\,(t_4-t_1^2)\ge2\,t_1^2\,\big(1-t_1\,t_3\,(1+t_2\,t_4)\big)>0\}
\end{multline*}
and
\begin{multline*}
\{u_1,u_2,u_3,u_4>0\mid
u_1>0,\;u_2>0,\;u_3>0,\;u_4>0,\;u_2^2\,u_3\,u_4<1,\\
u_2\,u_3\,(u_1-u_2)<1,\;u_1\,(1+u_2\,u_3)\le1\}.
\end{multline*}
\section{The two Wilson manifolds}
\label{49}
In this section I will present in detail the two charts for the two
Wilson manifolds, as introduced in \S\ref{25},
and corresponding to
the two graphs in Figure \ref{fig:3}.

\subsection{From Wilson to Hermite along Continuous Hahn}
Here we will see the chart for the Wilson manifold corresponding to the
first graph in Figure \ref{fig:3}.
The monic Wilson polynomial $w_n$ is given by \eqref{29}.
For $a_1,a_2,a_3,a_4>0$ put
\begin{align}
&p_n(x)=p_n(x;a_1,a_2,a_3,a_4):=
\rho^n\,w_n(\rho^{-1}x-\si;a,b,c,d),
\label{47}\sLP
{\rm where}\quad
&a=a_1^{-1}-\frac{1-a_1^{1/2}\,a_2\,a_4}{2\,a_1^{3/2}\,a_2^2\,a_3\,a_4}\,i,
\quad
b=a_1^{-1\,}a_2^{-1}
-\frac{1+a_1^{1/2}\,a_2\,a_4}{2\,a_1^{3/2}\,a_2^2\,a_3\,a_4}\,i,
\nonu\\
&c=a_1^{-1}+\frac{1-a_1^{1/2}\,a_2\,a_4}{2\,a_1^{3/2}\,a_2^2\,a_3\,a_4}\,i,
\quad
d=a_1^{-1}\,a_2{}^{-1}
+\frac{1+a_1^{1/2}\,a_2\,a_4}{2\,a_1^{3/2}\,a_2^2\,a_3\,a_4}\,i,
\nonu\\
&\rho=2^{3/2}\,a_1^2\,a_2^2\,a_3^2\,a_4,
\quad
\si=-\,\frac{1}{4\,a_1^3\,a_2^4\,a_3^2\,a_4^2}+
\frac{1-a_2}{2\,a_1^{5/2}\,a_2^3\,(1+a_2-a_1\,a_2)\,a_3^2\,a_4}\,.
\nonu
\end{align}
Then \eqref{32} holds with
\begin{align}
B_n&=2^{-1/2}
\Big(2\,n^4\,a_1^4\,a_2^4\,a_3^2\,a_4\,(1+a_2-a_1\,a_2)
+4\,n^3\,a_1^3\,a_2^3\,a_3^2\,a_4\,(2+2\,a_2-a_1\,a_2)\,(1+a_2-a_1\,a_2)
\nonu\\
&\qquad+n^2\,a_1^{3/2}\,a_2\,\big(2-2\,a_2+a_1^{1/2}\,a_2\,a_4\,
(1+a_2-a_1\,a_2)
+a_1^{1/2}\,a_2\,a_3^2\,a_4\,
(10+34\,a_2
\nonu\\
&\qquad\qquad-20\,a_1\,a_2+34\,a_2^2
-44\,a_1\,a_2^2+12\,a_1^2\,a_2^2+10\,a_2^3-20\,a_1\,a_2^3+12\,a_1^2\,a_2^3
-2\,a_1^3\,a_2^3)\big)
\nonu\displaybreak[0]\\
&\qquad+n\,a_1^{1/2}\,(2+2\,a_2-a_1\,a_2)\,
\big(2-2\,a_2+a_1^{1/2}\,a_2\,a_4\,
(1+ a_2-a_1\,a_2)
\nonu\\
&\qquad\qquad+2\,a_1^{1/2}\,a_2\,a_3^2\,a_4\,
(1+5\,a_2-2\,a_1\,a_2+5\,a_2^2
-6\,a_1\,a_2^2+a_1^2\,a_2^2+ a_2^3-2\,a_1\,a_2^3+a_1^2\,a_2^3)\big)
\nonu\\
&\qquad+(1+a_2-a_1\,a_2)\,\big(2\,a_1^{1/2}-2\,a_1^{1/2}\,a_2+a_4+2\,a_2\,a_4
-a_1\,a_2\,a_4+a_2^2\,a_4
\nonu\\
&\qquad\qquad-a_1\,a_2^2\,a_4
+4\,a_2\,a_3^2\,a_4\,(1+2 a_2-a_1\,a_2+a_2^2-a_1\,a_2^2)\big)\Big)
\nonu\\
&\quad\times\big(1+a_2-a_1\,a_2\big)^{-1}\,
\big(1+a_2+(n-1)\,a_1\,a_2\big)^{-1}\,
\big(1+a_2+n\,a_1\,a_2\big)^{-1},
\label{42}\mLP
C_n&=\thalf\,n\,
\Big(1+a_3^2\,\big(1+a_2+(n-1)\,a_1\,a_2\big)^2\Big)\,
\Big(1+a_1\,a_2^2\,a_3^2\,a_4^2\,\big(1+a_2+(n-1)\,a_1\,a_2\big)^2\Big)
\nonu\\
&\qquad\times\frac{(2-a_1+n\,a_1)\,(2+2\,a_2+(n-2)\,a_1\,a_2)\,
(2+(n-1)\,a_1\,a_2)}
{\big(1+a_2+(n-\thalf)\,a_1\,a_2\big)\,
\big(1+a_2+(n-1)\,a_1\,a_2\big)^2\,
\big(1+a_2+(n-\tfrac32)\,a_1\,a_2\big)}\,.
\label{43}
\end{align}
Note that $B_n$ and $C_n$, as functions of $a_1,a_2,a_3,a_4>0$,
can be uniquely extended to continuous functions of
$a_1,a_2,a_3,a_4\ge0$.

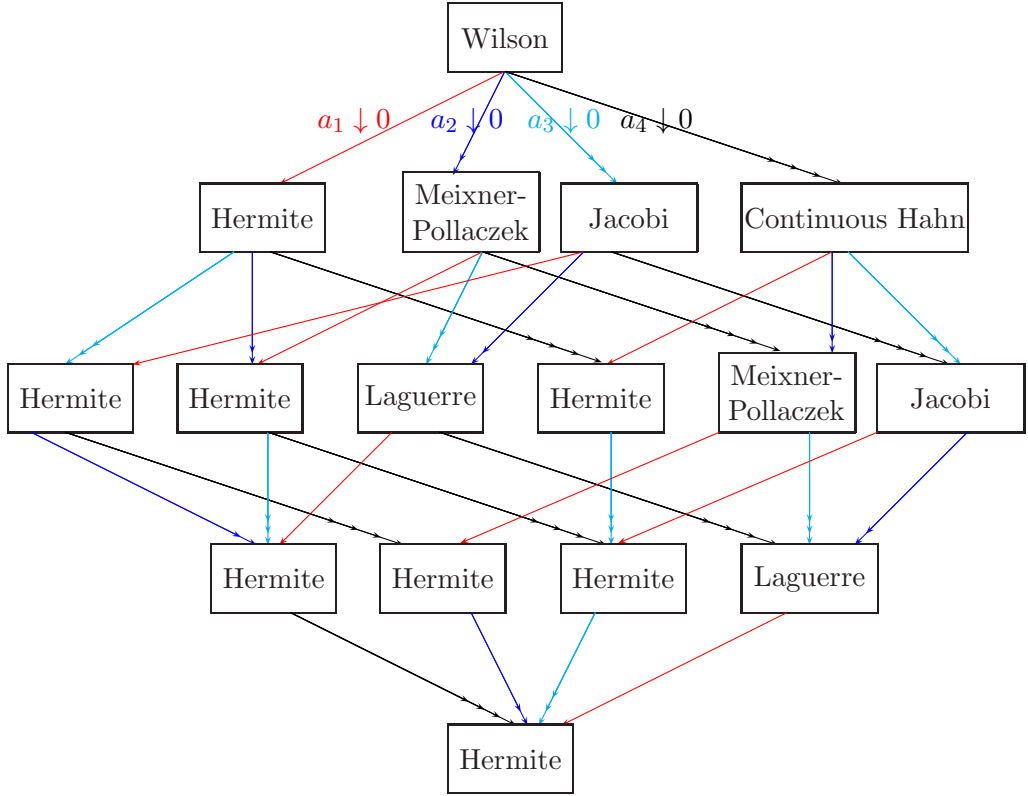
\begin{figure}[t]
\setlength{\unitlength}{3mm}
\begin{picture}(0,37.5)
\put(22,34.5) {\framebox(5,3) {Wilson}}
\put(24.5,34.5) {\color{red}\vector(-2,-1){10}}
\put(16.2,32){\color{red}{$a_1\downarrow0$}}
\put(24.5,34.5) {\color{blue}\vector(-1,-2){2.3}}
\put(24.5,34.5) {\color{blue}\vector(-1,-2){2.0}}
\put(21.2,32){\color{blue}{$a_2\downarrow0$}}
\put(24.5,34.5) {\color{cyan}\vector(1,-1){5}}
\put(24.5,34.5) {\color{cyan}\vector(1,-1){4.5}}
\put(24.5,34.5) {\color{cyan}\vector(1,-1){4}}
\put(25.5,32){\color{cyan}{$a_3\downarrow0$}}
\put(24.5,34.5) {\color{black}\vector(3,-1){15}}
\put(24.5,34.5) {\color{black}\vector(3,-1){14}}
\put(24.5,34.5) {\color{black}\vector(3,-1){13}}
\put(24.5,34.5) {\color{black}\vector(3,-1){12}}
\put(29.6,32){\color{black}{$a_4\downarrow0$}}
\put(11,26.5) {\framebox(5.5,3) {Hermite}}
\put(20,26.5) {\framebox(6,3.5) {\shortstack {Meixner-\\[1mm]Pollaczek}}}
\put(27,26.5) {\framebox(6,3) {Jacobi}}
\put(35,26.5) {\framebox(10,3) {Continuous Hahn}}
\put(12.5,26.5) {\color{cyan}\vector(-3,-2){7.5}}
\put(12.5,26.5) {\color{cyan}\vector(-3,-2){6.9}}
\put(12.5,26.5) {\color{cyan}\vector(-3,-2){6.3}}
\put(13.3,26.5) {\color{blue}\vector(0,-1){5}}
\put(13.3,26.5) {\color{blue}\vector(0,-1){4.5}}
\put(14.1,26.5) {\color{black}\vector(3,-1){14.7}}
\put(14.1,26.5) {\color{black}\vector(3,-1){13.7}}
\put(14.1,26.5) {\color{black}\vector(3,-1){12.7}}
\put(14.1,26.5) {\color{black}\vector(3,-1){11.7}}
\put(23.5,26.5) {\color{red}\vector(-2,-1){10}}
\put(23.5,26.5) {\color{cyan}\vector(-1,-2){2.5}}
\put(23.5,26.5) {\color{cyan}\vector(-1,-2){2.2}}
\put(23.5,26.5) {\color{cyan}\vector(-1,-2){1.9}}
\put(23.5,26.5) {\color{black}\vector(3,-1){13.2}}
\put(23.5,26.5) {\color{black}\vector(3,-1){12.2}}
\put(23.5,26.5) {\color{black}\vector(3,-1){11.2}}
\put(23.5,26.5) {\color{black}\vector(3,-1){10.2}}
\put(28,26.5) {\color{red}\vector(-4,-1){20}}
\put(28,26.5) {\color{blue}\vector(-1,-1){5}}
\put(28,26.5) {\color{blue}\vector(-1,-1){4.5}}
\put(29.2,26.5) {\color{black}\vector(3,-1){15}}
\put(29.2,26.5) {\color{black}\vector(3,-1){14}}
\put(29.2,26.5) {\color{black}\vector(3,-1){13}}
\put(29.2,26.5) {\color{black}\vector(3,-1){12}}
\put(39,26.5) {\color{red}\vector(-2,-1){10}}
\put(39,26.5) {\color{blue}\vector(0,-1){4.5}}
\put(39,26.5) {\color{blue}\vector(0,-1){4.0}}
\put(39.7,26.5) {\color{cyan}\vector(1,-1){5}}
\put(39.7,26.5) {\color{cyan}\vector(1,-1){4.5}}
\put(39.7,26.5) {\color{cyan}\vector(1,-1){4}}
\put(2.5,18.5) {\framebox(5.5,3) {Hermite}}
\put(10,18.5) {\framebox(5.5,3) {Hermite}}
\put(18,18.5) {\framebox(5.5,3) {Laguerre}}
\put(26,18.5) {\framebox(5.5,3) {Hermite}}
\put(34,18.5) {\framebox(6,3.5) {\shortstack {Meixner-\\[1mm]Pollaczek}}}
\put(41,18.5) {\framebox(6.5,3) {Jacobi}}
\put(3.5,18.5) {\color{blue}{\vector(2,-1){10}}}
\put(3.5,18.5) {\color{blue}{\vector(2,-1){9.3}}}
\put(5,18.5) {\color{black}{\vector(3,-1){15}}}
\put(5,18.5) {\color{black}{\vector(3,-1){14}}}
\put(5,18.5) {\color{black}{\vector(3,-1){13}}}
\put(5,18.5) {\color{black}{\vector(3,-1){12}}}
\put(14,18.5) {\color{cyan}\vector(0,-1){5}}
\put(14,18.5) {\color{cyan}\vector(0,-1){4.5}}
\put(14,18.5) {\color{cyan}\vector(0,-1){4}}
\put(14,18.5) {\color{black}\vector(3,-1){15}}
\put(14,18.5) {\color{black}\vector(3,-1){14}}
\put(14,18.5) {\color{black}\vector(3,-1){13}}
\put(14,18.5) {\color{black}\vector(3,-1){12}}
\put(19.5,18.5) {\color{red}\vector(-1,-1){5}}
\put(21.6,18.5) {\color{black}\vector(3,-1){15}}
\put(21.6,18.5) {\color{black}\vector(3,-1){14}}
\put(21.6,18.5) {\color{black}\vector(3,-1){13}}
\put(21.6,18.5) {\color{black}\vector(3,-1){12}}
\put(29.2,18.5) {\color{cyan}\vector(0,-1){5}}
\put(29.2,18.5) {\color{cyan}\vector(0,-1){4.5}}
\put(29.2,18.5) {\color{cyan}\vector(0,-1){4}}
\put(34,18.5) {\color{red}\vector(-7,-3){11.5}}
\put(38,18.5) {\color{cyan}\vector(0,-1){5}}
\put(38,18.5) {\color{cyan}\vector(0,-1){4.5}}
\put(38,18.5) {\color{cyan}\vector(0,-1){4}}
\put(41,18.5) {\color{red}\vector(-7,-3){11.5}}
\put(45,18.5) {\color{blue}\vector(-1,-1){5}}
\put(45,18.5) {\color{blue}\vector(-1,-1){4.5}}
\put(11.5,10.5) {\framebox(5.5,3) {Hermite}}
\put(19,10.5) {\framebox(5.5,3) {Hermite}}
\put(27,10.5) {\framebox(5.5,3) {Hermite}}
\put(35,10.5) {\framebox(6,3) {Laguerre}}
\put(15,10.5) {\color{black}\vector(2,-1){10}}
\put(15,10.5) {\color{black}\vector(2,-1){9.3}}
\put(15,10.5) {\color{black}\vector(2,-1){8.6}}
\put(15,10.5) {\color{black}\vector(2,-1){7.9}}
\put(23,10.5) {\color{blue}\vector(1,-2){2.5}}
\put(23,10.5) {\color{blue}\vector(1,-2){2.2}}
\put(28.5,10.5) {\color{cyan}\vector(-1,-2){2.5}}
\put(28.5,10.5) {\color{cyan}\vector(-1,-2){2.2}}
\put(28.5,10.5) {\color{cyan}\vector(-1,-2){1.9}}
\put(37,10.5) {\color{red}\vector(-2,-1){10}}
\put(22.0,2.5) {\framebox(5.5,3) {Hermite}}
\end{picture}
\vskip-1cm
\caption{The various parameter restrictions in the first Wilson manifold}
\label{fig:7}
\end{figure}

We now put one or more of the $a_1,a_2,a_3,a_4$ equal to zero
in \eqref{42} and \eqref{43} and
we proceed as in \S\ref{35}. We obtain:
\bLP
{\bf Continuous Hahn}:
\begin{align*}
&p_n(x;a_1,a_2,a_3,0)=
\rho^n\,p_n^{\rm CH}(\rho^{-1}x-\si;a,b,c,d),\sLP
&a=\frac{i+2\,a_2\,a_3}{2\,a_1\,a_2\,a_3}\,,\quad
b=\frac{-i+2 a_3}{2\,a_1\,a_2\,a_3}\,\quad
c=\frac{-i+2\,a_2\,a_3}{2\,a_1\,a_2\,a_3}\,,\quad
d=\frac{i+2\,a_3}{2\,a_1\,a_2\,a_3}\,,\\
&\rho=2^{3/2}\,a_1^{1/2}\,a_3,\quad
\si=\frac{1-a_2}{2\,a_1\,a_2\,a_3\,(1+a_2-a_1\,a_2)}\,.
\end{align*}
{\bf Jacobi}:
\begin{align*}
&p_n(x;a_1,a_2,0,a_4)=p_n(x-\si_0;a_1,a_2,0,0)=
\rho^n\,p_n^{(\al,\be)}(\rho^{-1}x-\si),\sLP
&\alpha=2a_1^{-1}-1,\quad
\beta=2\,a_1^{-1}\,a_2^{-1}-1,\quad
\rho=-2^{1/2}\,a_1^{-1/2}\,a_2^{-1},\\
&\si=-\,\frac{1-a_2}{1+a_2-a_1\,a_2}-\thalf\,a_1^{1/2}\,a_2\,a_4\,,\quad
\si_0=2^{-1/2}\,a_4.
\end{align*}
{\bf Meixner-Pollaczek}:
\begin{align*}
&p_n(x;a_1,0,a_3,a_4)=p_n(x-\si_0;a_1,0,a_3,0)=
\rho^n\,p_n^{(\la)}(\rho^{-1}x-\si),\sLP
&\lambda=a_1{}^{-1},\quad
\phi=\arctan(a_3),\quad
\rho=-2^{3/2}\,a_1^{1/2}\,a_3,\quad
\si=\frac{2-a_1}{2\,a_1\,a_3}-\tfrac14\,a_1^{-1/2}\,a_3^{-1}\,a_4\,,\quad
\si_0=2^{-1/2}\,a_4.
\end{align*}
{\bf Laguerre}:
\begin{align*}
&p_n(x;a_1,0,0,a_4)=p_n(x-\si_0;a_1,0,0,0)=
\rho^n\,\ell_n^{(\al)}(\rho^{-1}x-\si),\sLP
&\al=2a_1^{-1}-1,\quad
\rho=2^{1/2}\,a_1^{1/2},\quad
\si=1-2a_1^{-1}+\thalf\,a_1^{-1/2}\,a_4,\quad
\si_0=2^{-1/2}\,a_4.
\end{align*}
{\bf Hermite}:
\begin{align*}
&p_n(x;0,a_2,a_3,a_4)=p_n(x-\si_0;0,a_2,a_3,0)
=\rho_1^n\,p_n(\rho_1^{-1}\,x-\si_1;0,a_2,0,a_4)\\
&\hskip2.9cm=\rho_2^n\,p_n(\rho_2^{-1}\,x-\si_2;0,0,a_3,a_4)
=\rho_1^n\,p_n(\rho_1^{-1}\,x-\si_3;0,a_2,0,0)\\
&\hskip2.9cm=\rho^n\,p_n(\rho^{-1}\,x-\si_4;0,0,a_3,0)
=\rho_2^n\,p_n(\rho_2^{-1}\,x-\si_5;0,0,0,a_4)\\
&\hskip2.9cm=(2^{-3/2}\rho)^n\,p_n\big(2^{3/2}(\rho^{-1}\,x-\si);0,0,0,0\big)
=\rho^n\,h_n(\rho^{-1}x-\si),\sLP
&\rho=\frac{2^{3/2}\,(1+(1+a_2)^2\,a_3^2)^{1/2}}{(1+a_2)^{3/2}}\,,\quad
\si=\frac{a_4\,(1+a_2)^{3/2}\,(1+4\,a_2\,a_3^2)}
{4\,(1+(1+a_2)^2\,a_3^2)^{1/2}}\,,\quad
\si_0=2^{-1/2}\,(1+4\,a_2\,a_3^2)\,a_4,\\
&\rho_1=(1+(1+a_2)^2\,a_3^2)^{1/2},\quad
\si_1=\frac{a_4\,\big(1+4\,a_2\,a_3^2-
(1+(1+a_2)^2\,a_3^2)^{1/2}\big)}
{2^{1/2}\,(1+(1+a_2)^2\,a_3^2)^{1/2}}\,,\\
&\rho_2=\frac{(1+(1+a_2)^2\,a_3^2)^{1/2}}{(1+a_2)^{3/2}\,
(1+a_3^2)^{1/2}}\,,\quad
\si_2=\frac{a_4\,(1+a_2)^{3/2}\,(1+a_3^2)^{1/2}\,(1+4\,a_2\,a_3^2)}
{2^{1/2}\,(1+(1+a_2)^2\,a_3^2)^{1/2}}\,-2^{-1/2}\,a_4,\\
&\si_3=\frac{a_4\,(1+4\,a_2\,a_3^2)}
{2^{1/2}\,(1+(1+a_2)^2\,a_3^2)^{1/2}}\,,\quad
\si_4=\frac{a_4\,(1+a_2)^{3/2}\,(1+a_3^2)^{1/2}\,(1+4\,a_2\,a_3^2)}
{2^{1/2}\,(1+(1+a_2)^2\,a_3^2)^{1/2}}\,,\\
&\si_5=\frac{a_4\,(1+a_2)^{3/2}\,(1+4\,a_2\,a_3^2)}
{2^{1/2}\,(1+(1+a_2)^2\,a_3^2)^{1/2}}\,-2^{-1/2}\,a_4.
\end{align*}
The various parameter restrictions of the polynomial \eqref{47}
are summarized in Figure \ref{fig:7}. Many boxes essentially coincide.
\subsection{From Wilson to Hermite along Continuous Dual Hahn}
Here we will see the chart for the Wilson manifold corresponding to the
first graph in Figure \ref{fig:3}.
The monic Wilson polynomial $w_n$ is given by \eqref{29}.
For $b_1,b_2,b_3,b_4>0$ put
\begin{align}
&p_n(x)=p_n(x;b_1,b_2,b_3,b_4):=
\rho^n\,w_n(\rho^{-1}x-\si;a,b,c,d),
\label{44}\\
\intertext{where}
&a=\frac{1+b_1}{2\,b_1}+\frac{1+4\,b_1\,b_2}{2\,b_1^3\,b_2\,b_3}\,i,
\quad
b=\frac{1+b_1}{2\,b_1}-\frac{1+4\,b_1\,b_2}{2\,b_1^3\,b_2\,b_3}\,i,\quad
c=\frac{1}{b_1^6\,b_2^2\,b_3^3\,b_4}
\nonu\\
&+\frac{2+b_1\,b_3+b_1^3\,b_2\,b_3^2
+3\,b_1^4\,b_2\,b_3^2}{2\,b_1^4\,b_2\,b_3^2}\,,
\quad
d=-\frac{2+b_1\,b_3+b_1^3\,b_2\,b_3^2+b_1^4\,b_2\,b_3^2}{2\,b_1^4\,b_2\,b_3^2}
\,,\quad
\rho=2^{-1/2}\,b_1^{9/2}\,b_2\,b_3^2,
\nonu\\
&\si=-\,\frac{1+4\,b_1\,b_2+4\,b_1^2\,b_2\,b_4+4\,b_1^3\,b_2^2\,b_3\,b_4
+4\,b_1^4\,b_2^2\,b_3^2\,b_4+4\,b_1^4\,b_2^2\,b_3\,b_4^2
+4\,b_1^5\,b_2^2\,b_3^2\,b_4^2}
{4\,b_1^6\,b_2^2\,b_3^2}\,.\nonu
\end{align}
Then \eqref{32} holds with
\begin{align}
&B_n=
2^{-5/2}b_1^{1/2}\,
\Big(8\,n^4\,b_1^{16}\,b_2^5\,b_3^8\,b_4^2
+16\,n^3\,b_1^{10}\,b_2^3\,b_3^5\,b_4\,
(1+b_1^5\,b_2^2\,b_3^3\,b_4+b_1^6\,b_2^2\,b_3^3\,b_4)
\nonu\displaybreak[0]\\
&\qquad
+4\,n^2\,b_1^4\,b_2\,b_3^2\,
\big(2-2\,b_1^2\,b_2\,b_3\,b_4-b_1^3\,b_2\,b_3^2\,b_4
+5\,b_1^5\,b_2^2\,b_3^3\,b_4
+4\,b_1^6\,b_2^2\,b_3^3\,b_4-2\,b_1^4\,b_2^2\,b_3^2\,b_4^2
\nonu\displaybreak[0]\\
&\qquad\qquad-2\,b_1^5\,b_2^2\,b_3^3\,b_4^2-b_1^6\,b_2^2\,b_3^4\,b_4^2
-2\,b_1^7\,b_2^3\,b_3^4\,b_4^2-4\,b_1^8\,b_2^3\,b_3^4\,b_4^2
+8\,b_1^8\,b_2^4\,b_3^4\,b_4^2
-b_1^8\,b_2^3\,b_3^5\,b_4^2\nonu\\
&\qquad\qquad
-2\,b_1^9\,b_2^3\,b_3^5\,b_4^2
+2\,b_1^{10}\,b_2^4\,b_3^6\,b_4^2+4\,b_1^{11}\,b_2^4\,b_3^6\,b_4^2
+b_1^{12}\,b_2^4\,b_3^6\,b_4^2
-4\,b_1^8\,b_2^3\,b_3^4\,b_4^3-4\,b_1^9\,b_2^4\,b_3^5\,b_4^3\nonu\\
&\qquad\qquad
-4\,b_1^{10}\,b_2^4\,b_3^6\,b_4^3-4\,b_1^{10}\,b_2^4\,b_3^5\,b_4^4
-4\,b_1^{11}\,b_2^4\,b_3^6\,b_4^4\big)
\nonu\displaybreak[0]\\
&\qquad
-4\,n\,(1+b_1^5\,b_2^2\,b_3^3\,b_4+b_1^6\,b_2^2\,b_3^3\,b_4)\,
\big(2+b_1\,b_3-b_1^3\,b_2\,b_3^2+2\,b_1^2\,b_2\,b_3\,b_4
+2\,b_1^3\,b_2\,b_3^2\,b_4\nonu\\
&\qquad\qquad+b_1^4\,b_2\,b_3^3\,b_4+2\,b_1^5\,b_2^2\,b_3^3\,b_4
+4\,b_1^6\,b_2^2\,b_3^3\,b_4-8\,b_1^6\,b_2^3\,b_3^3\,b_4
+b_1^6\,b_2^2\,b_3^4\,b_4+2\,b_1^7\,b_2^2\,b_3^4\,b_4\nonu\\
&\qquad\qquad
+b_1^{10}\,b_2^3\,b_3^5\,b_4+4\,b_1^6\,b_2^2\,b_3^3\,b_4^2
+4\,b_1^7\,b_2^3\,b_3^4\,b_4^2+4\,b_1^8\,b_2^3\,b_3^5\,b_4^2
+4\,b_1^8\,b_2^3\,b_3^4\,b_4^3+4\,b_1^9\,b_2^3\,b_3^5\,b_4^3\big)\nonu
\displaybreak[0]\\
&\qquad-\big(1+b_1^5\,b_2^2\,b_3^3\,b_4)
(4-16\,b_2+2\,b_3+2\,b_1\,b_3+b_1^2\,b_2\,b_3^2+2\,b_1^3\,b_2\,b_3^2
+b_1^4\,b_2\,b_3^2+4\,b_4\nonu\\
&\qquad\qquad
+8\,b_1\,b_2\,b_3\,b_4+4\,b_1^2\,b_2\,b_3\,b_4+8\,b_1^2\,b_2\,b_3^2\,b_4
+4\,b_1^3\,b_2\,b_3^2\,b_4+2\,b_1^3\,b_2\,b_3^3\,b_4+2\,b_1^4\,b_2\,b_3^3\,b_4
\nonu\\
&\qquad\qquad
+8\,b_1^4\,b_2^2\,b_3^3\,b_4+12\,b_1^5\,b_2^2\,b_3^3\,b_4
+8\,b_1^6\,b_2^2\,b_3^3\,b_4-16\,b_1^6\,b_2^3\,b_3^3\,b_4
+2\,b_1^5\,b_2^2\,b_3^4\,b_4+6\,b_1^6\,b_2^2\,b_3^4\,b_4\nonu\\
&\qquad\qquad
+4\,b_1^7\,b_2^2\,b_3^4\,b_4+b_1^7\,b_2^3\,b_3^5\,b_4
+4\,b_1^8\,b_2^3\,b_3^5\,b_4
+5\,b_1^9\,b_2^3\,b_3^5\,b_4+2\,b_1^{10}\,b_2^3\,b_3^5\,b_4
+4\,b_1^2\,b_2\,b_3\,b_4^2\nonu\\
&\qquad\qquad
+4\,b_1^3\,b_2\,b_3^2\,b_4^2
+4\,b_1^5\,b_2^2\,b_3^3\,b_4^2+8\,b_1^6\,b_2^2\,b_3^3\,b_4^2
+4\,b_1^6\,b_2^3\,b_3^4\,b_4^2
+8\,b_1^7\,b_2^3\,b_3^4\,b_4^2
+4\,b_1^7\,b_2^3\,b_3^5\,b_4^2\nonu\\
&\qquad\qquad+8\,b_1^8\,b_2^3\,b_3^5\,b_4^2
+4\,b_1^7\,b_2^3\,b_3^4\,b_4^3
+8\,b_1^8\,b_2^3\,b_3^4\,b_4^3
+4\,b_1^8\,b_2^3\,b_3^5\,b_4^3+8\,b_1^9\,b_2^3\,b_3^5\,b_4^3\big)\Big)
\nonu\displaybreak[0]\\
&\qquad\times
\big(1+b_1^5\,b_2^2\,b_3^3\,b_4+2\,n\,b_1^6\,b_2^2\,b_3^3\,b_4\big)^{-1}
\big(1+b_1^5\,b_2^2\,b_3^3\,b_4+2\,(1+n)\,b_1^6\,b_2^2\,b_3^3\,b_4\big)^{-1},
\label{45}
\displaybreak[0]\mLP
&C_n=\frac{1}{8}\,n\,(1+n\,b_1)\,(1+n\,b_1^6\,b_2^2\,b_3^3\,b_4)\,
\big(1+b_1^5\,b_2^2\,b_3^3\,b_4\,(1+n\,b_1)\big)
\nonu\\
&\times
\frac{2+2\,b_1\,b_3+b_1^2\,b_3^2+4\,b_1^3\,b_2\,b_3^2
+8\,b_1^4\,b_2^2\,b_3^2
+2\,(1-n)\,b_1^4\,b_2\,b_3^2\,(2+b_1\,b_3)+2\,(1-n)^2\,b_1^8\,b_2^2\,b_3^4}
{\big(1+b_1^5\,b_2^2\,b_3^3\,b_4\,(1+(2n-1)\,b_1)\big)\,
\big(1+b_1^5\,b_2^2\,b_3^3\,b_4\,(1+2\,n\,b_1)\big)^2\,
\big(1+b_1^5\,b_2^2\,b_3^3\,b_4\,(1+(2\,n+1)\,b_1)\big)}\nonu\\
&
\times\big(2+4\,b_1^2\,b_2\,b_3\,b_4+2\,b_1^3\,b_2\,b_3^2\,b_4
+4\,b_1^5\,b_2^2\,b_3^3\,b_4+4\,(1+n)\,b_1^6\,b_2^2\,b_3^3\,b_4
+2\,b_1^4\,b_2^2\,b_3^2\,b_4^2\nonu\\
&\qquad\qquad+2\,b_1^5\,b_2^2\,b_3^3\,b_4^2+b_1^6\,b_2^2\,b_3^4\,b_4^2
+8\,b_1^7\,b_2^3\,b_3^4\,b_4^2+4\,(1+n)\,b_1^8\,b_2^3\,b_3^4\,b_4^2
+8\,b_1^8\,b_2^4\,b_3^4\,b_4^2\nonu\\
&\qquad\qquad+2b_1^8\,b_2^3\,b_3^5\,b_4^2\,(1+(n+1)\,b_1)
+2\,b_1^{10}\,b_2^4\,b_3^6\,b_4^2\,(1+(n+1)\,b_1)^2\big).
\label{46}
\end{align}
Note that $B_n$ and $C_n$, as functions of $b_1,b_2,b_3,b_4>0$,
can be uniquely extended to continuous functions of
$b_1,b_2,b_3,b_4\ge0$.

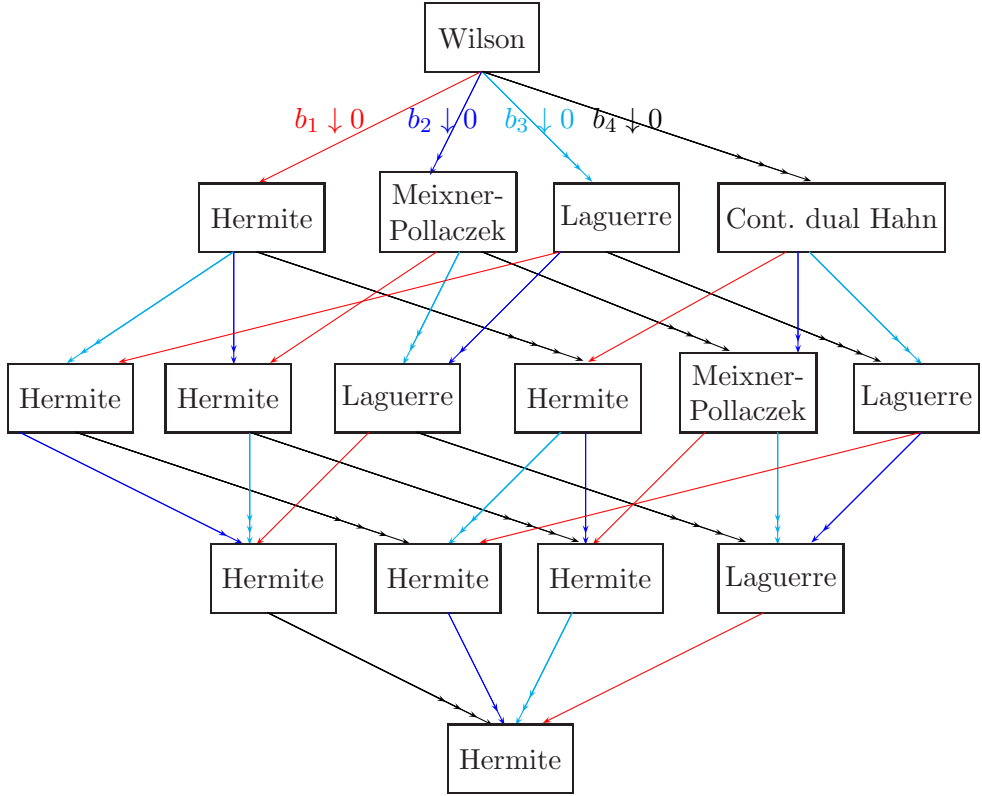
\begin{figure}[t]
\setlength{\unitlength}{3mm}
\begin{picture}(0,37.5)
\put(22,34.5) {\framebox(5,3) {Wilson}}
\put(24.5,34.5) {\color{red}\vector(-2,-1){9.9}}
\put(16.2,32){\color{red}{$b_1\downarrow0$}}
\put(24.5,34.5) {\color{blue}\vector(-1,-2){2.3}}
\put(24.5,34.5) {\color{blue}\vector(-1,-2){2.0}}
\put(21.2,32){\color{blue}{$b_2\downarrow0$}}
\put(24.5,34.5) {\color{cyan}\vector(1,-1){4.9}}
\put(24.5,34.5) {\color{cyan}\vector(1,-1){4.4}}
\put(24.5,34.5) {\color{cyan}\vector(1,-1){3.9}}
\put(25.5,32){\color{cyan}{$b_3\downarrow0$}}
\put(24.5,34.5) {\color{black}\vector(3,-1){14.6}}
\put(24.5,34.5) {\color{black}\vector(3,-1){13.6}}
\put(24.5,34.5) {\color{black}\vector(3,-1){12.6}}
\put(24.5,34.5) {\color{black}\vector(3,-1){11.6}}
\put(29.4,32){\color{black}{$b_4\downarrow0$}}
\put(12,26.5) {\framebox(5.5,3) {Hermite}}
\put(20,26.5) {\framebox(6,3.5) {\shortstack {Meixner-\\[1mm]Pollaczek}}}
\put(27.7,26.5) {\framebox(5.5,3) {Laguerre}}
\put(35,26.5) {\framebox(10,3) {Cont.\ dual Hahn}}
\put(13.5,26.5) {\color{cyan}\vector(-3,-2){7.4}}
\put(13.5,26.5) {\color{cyan}\vector(-3,-2){6.8}}
\put(13.5,26.5) {\color{cyan}\vector(-3,-2){6.2}}
\put(13.5,26.5) {\color{blue}\vector(0,-1){5}}
\put(13.5,26.5) {\color{blue}\vector(0,-1){4.5}}
\put(14.5,26.5) {\color{black}\vector(3,-1){14.5}}
\put(14.5,26.5) {\color{black}\vector(3,-1){13.5}}
\put(14.5,26.5) {\color{black}\vector(3,-1){12.5}}
\put(14.5,26.5) {\color{black}\vector(3,-1){11.5}}
\put(22.5,26.5) {\color{red}\vector(-3,-2){7.4}}
\put(23.5,26.5) {\color{cyan}\vector(-1,-2){2.5}}
\put(23.5,26.5) {\color{cyan}\vector(-1,-2){2.2}}
\put(23.5,26.5) {\color{cyan}\vector(-1,-2){1.9}}
\put(24.5,26.5) {\color{black}\vector(5,-2){11}}
\put(24.5,26.5) {\color{black}\vector(5,-2){10.2}}
\put(24.5,26.5) {\color{black}\vector(5,-2){9.4}}
\put(24.5,26.5) {\color{black}\vector(5,-2){8.6}}
\put(28,26.5) {\color{red}\vector(-4,-1){19.6}}
\put(28,26.5) {\color{blue}\vector(-1,-1){5}}
\put(28,26.5) {\color{blue}\vector(-1,-1){4.5}}
\put(30,26.5) {\color{black}\vector(5,-2){12.2}}
\put(30,26.5) {\color{black}\vector(5,-2){11.4}}
\put(30,26.5) {\color{black}\vector(5,-2){10.6}}
\put(30,26.5) {\color{black}\vector(5,-2){9.8}}
\put(38,26.5) {\color{red}\vector(-9,-5){8.8}}
\put(38.5,26.5) {\color{blue}\vector(0,-1){4.5}}
\put(38.5,26.5) {\color{blue}\vector(0,-1){4.0}}
\put(39,26.5) {\color{cyan}\vector(1,-1){5}}
\put(39,26.5) {\color{cyan}\vector(1,-1){4.5}}
\put(39,26.5) {\color{cyan}\vector(1,-1){4}}
\put(3.5,18.5) {\framebox(5.5,3) {Hermite}}
\put(10.5,18.5) {\framebox(5.5,3) {Hermite}}
\put(18,18.5) {\framebox(5.5,3) {Laguerre}}
\put(26,18.5) {\framebox(5.5,3) {Hermite}}
\put(33.3,18.5) {\framebox(6,3.5) {\shortstack {Meixner-\\[1mm]Pollaczek}}}
\put(41,18.5) {\framebox(5.5,3) {Laguerre}}
\put(4,18.5) {\color{blue}{\vector(2,-1){9.9}}}
\put(4,18.5) {\color{blue}{\vector(2,-1){9.2}}}
\put(6.5,18.5) {\color{black}{\vector(3,-1){14.8}}}
\put(6.5,18.5) {\color{black}{\vector(3,-1){13.8}}}
\put(6.5,18.5) {\color{black}{\vector(3,-1){12.8}}}
\put(6.5,18.5) {\color{black}{\vector(3,-1){11.8}}}
\put(14.2,18.5) {\color{cyan}\vector(0,-1){5}}
\put(14.2,18.5) {\color{cyan}\vector(0,-1){4.5}}
\put(14.2,18.5) {\color{cyan}\vector(0,-1){4}}
\put(14.2,18.5) {\color{black}\vector(3,-1){14.6}}
\put(14.2,18.5) {\color{black}\vector(3,-1){13.6}}
\put(14.2,18.5) {\color{black}\vector(3,-1){12.6}}
\put(14.2,18.5) {\color{black}\vector(3,-1){11.6}}
\put(19.5,18.5) {\color{red}\vector(-1,-1){5}}
\put(21.6,18.5) {\color{black}\vector(3,-1){14.6}}
\put(21.6,18.5) {\color{black}\vector(3,-1){13.6}}
\put(21.6,18.5) {\color{black}\vector(3,-1){12.6}}
\put(21.6,18.5) {\color{black}\vector(3,-1){11.6}}
\put(28,18.5) {\color{cyan}\vector(-1,-1){5}}
\put(28,18.5) {\color{cyan}\vector(-1,-1){4.5}}
\put(28,18.5) {\color{cyan}\vector(-1,-1){4}}
\put(29.1,18.5) {\color{blue}\vector(0,-1){5}}
\put(29.1,18.5) {\color{blue}\vector(0,-1){4.5}}
\put(34.4,18.5) {\color{red}\vector(-1,-1){5}}
\put(37.6,18.5) {\color{cyan}\vector(0,-1){5}}
\put(37.6,18.5) {\color{cyan}\vector(0,-1){4.5}}
\put(37.6,18.5) {\color{cyan}\vector(0,-1){4}}
\put(44,18.5) {\color{red}\vector(-4,-1){19.6}}
\put(44,18.5) {\color{blue}\vector(-1,-1){4.9}}
\put(44,18.5) {\color{blue}\vector(-1,-1){4.4}}
\put(12.5,10.5) {\framebox(5.5,3) {Hermite}}
\put(19.8,10.5) {\framebox(5.5,3) {Hermite}}
\put(27,10.5) {\framebox(5.5,3) {Hermite}}
\put(35,10.5) {\framebox(5.5,3) {Laguerre}}
\put(15,10.5) {\color{black}\vector(2,-1){10}}
\put(15,10.5) {\color{black}\vector(2,-1){9.3}}
\put(15,10.5) {\color{black}\vector(2,-1){8.6}}
\put(15,10.5) {\color{black}\vector(2,-1){7.9}}
\put(23,10.5) {\color{blue}\vector(1,-2){2.5}}
\put(23,10.5) {\color{blue}\vector(1,-2){2.2}}
\put(28.5,10.5) {\color{cyan}\vector(-1,-2){2.5}}
\put(28.5,10.5) {\color{cyan}\vector(-1,-2){2.2}}
\put(28.5,10.5) {\color{cyan}\vector(-1,-2){1.9}}
\put(37,10.5) {\color{red}\vector(-2,-1){9.8}}
\put(23,2.5) {\framebox(5.5,3) {Hermite}}
\end{picture}
\vskip-1cm
\caption{The various parameter restrictions in the second Wilson manifold}
\label{fig:8}
\end{figure}

We now put one or more of the $b_1,b_2,b_3,b_4$ equal to zero
in \eqref{45} and \eqref{46} and
we proceed as in \S\ref{35}. We obtain:
\bLP
{\bf Continuous dual Hahn}:
\begin{align*}
&p_n(x;b_1,b_2,b_3,0)=\rho^n\,s_n(\rho^{-1}x-\si;a,b,c),\\
&a=\frac{1+b_1}{2\,b_1}+\frac{1+4\,b_1\,b_2}{2\,b_1^3\,b_2\,b_3}\,i,\quad
b=\frac{1+b_1}{2\,b_1}-\frac{1+4\,b_1\,b_2}{2b_1^3\,b_2\,b_3}\,i,\\
&c=-\frac{2+b_1\,b_3+b_1^3\,b_2\,b_3^2+b_1^4\,b_2\,b_3^2}{2\,b_1^4\,b_2\,b_3^2},
\quad
\rho=2^{-1/2}\,b_1^{9/2}\,b_2\,b_3^2,\quad
\si=-\frac{1+4\,b_1\,b_2}{4\,b_1^6\,b_2^2\,b_3^2}\,.
\end{align*}
{\bf Meixner-Pollaczek}:
\begin{align*}
&p_n(x;b_1,0,b_3,b_4)=p_n(x-\si_0;b_1,0,b_3,0)=
\rho^n\,p_n^{(\la)}(\rho^{-1}x-\si),\sLP
&\la=\frac{1+b_1}{2\,b_1}\,,\quad
\phi=\arctan\Big(\frac{b_1\,b_3}{2+b_1\,b_3}\Big),\quad
\rho=2^{-1/2}\,b_1^{3/2}\,b_3\,\quad
\si=\frac{1-b_1\,b_4}{b_1^2\,b_3}\,,\quad
\si_0=-2^{-1/2}\,b_1^{1/2}\,b_4.
\end{align*}
{\bf Laguerre}:
\begin{align*}
&p_n(x;b_1,b_2,0,b_4)=p_n(x-\si_0;b_1,0,0,b_4)=
p_n(x-\si_1;b_1,b_2,0,0)\\
&=p_n(x-\si_1;b_1,0,0,0)
=\rho^n\,\ell_n^{(\al)}(\rho^{-1}x-\si),\sLP
&\al=b_1^{-1},\quad
\rho=2^{-1/2}\,b_1^{1/2},\quad
\si=b_4-4b_2-b_1^{-1},\quad
\si_0=2^{3/2}\,b_1^{1/2}\,b_2,\quad
\si_1=-2^{-1/2}\,b_1^{1/2}\,b_4.
\end{align*}
{\bf Hermite}:
\begin{align*}
&p_n(x;0,b_2,b_3,b_4)=p_n(x;0,0,0,0)=\rho^n\,h_n(\rho^{-1}x-\si),\sLP
&\rho=1,\quad\si=0.
\end{align*}
The various parameter restrictions of the polynomial \eqref{44}
are summarized in Figure \ref{fig:8}. Many boxes almost coincide.
\section{Discussion of the results}
Various aspects of the results presented in this paper need further
attention in future. I list them below.
\begin{enumerate}
\item
The various reparametrizations (three for the Racah polynomials and
two for the Wilson polynomials) were obtained by computer algebra
experiments: by varying possible reparametrizations until they ``fit'', i.e,
until the coefficients $B_n$ and $C_n$ in the recurrence relation \eqref{32}
are continuous on all parts of the boundary and give there
$B_n$ and $C_n$ for families on lower levels of the Askey scheme.
Is each of these reparametrizations canonical in some way, in the sense
that it is unique, or generates all other reparametrizations which
do the same job?
\item
Do we really need three charts for covering the Racah manifold?
In particular, does there exist a ``fusion'' of the second and third
chart for the Racah manifold such that everything below dual Hahn,
including Laguerre, is covered?
\item
Gluing between the two Wilson manifolds should be considered on
the three-dimensional part corresponding with the three-parameter
family of monic Wilson polynomials $w_n(x;a,b,c,d)$ having one pair
of complex conjugate parameters and the other two parameters equal and real.
Also, gluing on lower dimensional parts between Wilson manifolds and
the Racah manifold should be considered such that one or more of the
boxes for Jacobi, Laguerre and Hermite in the various manifolds are
identified.
\item
As we have seen in Figures \ref{fig:4}--\ref{fig:8}, some families occur
in various parts of the boundary in a given chart. Hence, our geometric
description of the Askey scheme is not in bijective correspondence with
the Askey scheme, but this will only be achieved by considering suitable
quotient spaces of our four-manifolds with corners. A precise description
of these quotient spaces should be worked out.
\item
An obvious challenge would be to extend the results of this paper to
the $q$-Askey scheme, including the limits for $q\uparrow1$. A trial for a
small part of the $q$-Askey scheme with its limits to $q=1$ would already
be interesting.
\end{enumerate}

\quad\\
\begin{footnotesize}
\begin{quote}
{T. H. Koornwinder, Korteweg-de Vries Institute, University of
Amsterdam,\mLP
P.O.\ Box 94248, 1090 GE Amsterdam, The Netherlands;
\mLP
email: }{\tt T.H.Koornwinder@uva.nl}
\end{quote}
\end{footnotesize}

\end{document}